\documentclass[11pt]{article}

\usepackage{amsfonts,amssymb,amsmath,amsthm,latexsym}
\usepackage{mathrsfs}
\usepackage{graphicx,graphics,psfrag,epsf}

\usepackage[english]{babel}
 \makeatletter
 \@addtoreset{equation}{section}
 \makeatother

 \makeatletter
 \@addtoreset{enunciato}{section}
 \makeatother

 \newcounter{enunciato}[section]

 \newtheorem{ittheorem}{Theorem}
 \newtheorem{itlemma}{Lemma}
 \newtheorem{itproposition}{Proposition}
 \newtheorem{itdefinition}{Definition}
 \newtheorem{itremark}{Remark}
 \newtheorem{itclaim}{Claim}
 \newtheorem{itfact}{Fact}
 \newtheorem{itconjecture}{Conjecture}
 \newtheorem{itcorollary}{Corollary}

 \newenvironment{theorem}{\addtocounter{enunciato}{1}
 \begin{ittheorem}}{\end{ittheorem}}

 \newenvironment{lemma}{\addtocounter{enunciato}{1}
 \begin{itlemma}}{\end{itlemma}}

 \newenvironment{proposition}{\addtocounter{enunciato}{1}
 \begin{itproposition}}{\end{itproposition}}

 \newenvironment{definition}{\addtocounter{enunciato}{1}
 \begin{itdefinition}}{\end{itdefinition}}

 \newenvironment{remark}{\addtocounter{enunciato}{1}
 \begin{itremark}}{\end{itremark}}

 \newenvironment{claim}{\addtocounter{enunciato}{1}
 \begin{itclaim}}{\end{itclaim}}

 \newenvironment{fact}{\addtocounter{enunciato}{1}
 \begin{itfact}}{\end{itfact}}

 \newenvironment{conjecture}{\addtocounter{enunciato}{1}
 \begin{itconjecture}}{\end{itconjecture}}

 \newenvironment{corollary}{\addtocounter{enunciato}{1}
 \begin{itcorollary}}{\end{itcorollary}}

 \newcommand{\be}[1]{\begin{equation}\label{#1}}
 \newcommand{\ee}{\end{equation}}

 \newcommand{\bl}[1]{\begin{lemma}\label{#1}}
 \newcommand{\el}{\end{lemma}}

 \newcommand{\br}[1]{\begin{remark}\label{#1}}
 \newcommand{\er}{\end{remark}}

 \newcommand{\bt}[1]{\begin{theorem}\label{#1}}
 \newcommand{\et}{\end{theorem}}

 \newcommand{\bd}[1]{\begin{definition}\label{#1}}
 \newcommand{\ed}{\end{definition}}

 \newcommand{\bcl}[1]{\begin{claim}\label{#1}}
 \newcommand{\ecl}{\end{claim}}

 \newcommand{\bfact}[1]{\begin{fact}\label{#1}}
 \newcommand{\efact}{\end{fact}}

 \newcommand{\bp}[1]{\begin{proposition}\label{#1}}
 \newcommand{\ep}{\end{proposition}}

 \newcommand{\bc}[1]{\begin{corollary}\label{#1}}
 \newcommand{\ec}{\end{corollary}}

 \newcommand{\bcj}[1]{\begin{conjecture}\label{#1}}
 \newcommand{\ecj}{\end{conjecture}}

 \newcommand{\bpr}{\begin{proof}}
 \newcommand{\epr}{\end{proof}}

 \newcommand{\bprl}[1]{\begin{proofof}{\it\ref{#1}}.\,\,}
 \newcommand{\eprl}{\end{proofof}}

 \newcommand{\bi}{\begin{itemize}}
 \newcommand{\ei}{\end{itemize}}

 \newcommand{\ben}{\begin{enumerate}}
 \newcommand{\een}{\end{enumerate}}

\def \Z {\mathbb{Z}}
\def \N {\mathbb{N}}

\begin{document}


\title{Law of large numbers for non-elliptic\\ 
random walks in dynamic random environments}

\author{\renewcommand{\thefootnote}{\arabic{footnote}}
F.\ den Hollander \footnotemark[1]\,\,\,\,\footnotemark[2]
\\
\renewcommand{\thefootnote}{\arabic{footnote}}
R.\ dos Santos \footnotemark[1]\,\,\,\,\footnotemark[5]
\\
\renewcommand{\thefootnote}{\arabic{footnote}}
V.\ Sidoravicius \footnotemark[3]\,\,\,\,\footnotemark[4]}

\footnotetext[1]{
Mathematical Institute, Leiden University, P.O.\ Box 9512,
2300 RA Leiden, The Netherlands.}

\footnotetext[2]{
EURANDOM, P.O.\ Box 513, 5600 MB Eindhoven, The Netherlands.}

\footnotetext[3]{
CWI, Science Park 123, 1098 XG, Amsterdam, The Netherlands.}

\footnotetext[4]{
IMPA, Estrada Dona Castorina 110, Jardim Botanico, 
CEP 22460-320, Rio de Janeiro, Brasil.}

\footnotetext[5]{Corresponding author. Email: renato@math.leidenuniv.nl, phone: +31 71 527-7141, fax: +31 71 527-7101.}

\maketitle

\begin{abstract}
We prove a law of large numbers for a class of $\Z^d$-valued random walks 
in dynamic random environments, including non-elliptic examples. 
We assume for the random environment a mixing property called \emph{conditional cone-mixing} 
and that the random walk tends to stay inside wide enough space-time cones. 
The proof is based on a generalization of a regeneration scheme developed by 
Comets and Zeitouni~\cite{CoZe} for static random environments and adapted by 
Avena, den Hollander and Redig~\cite{AvdHoRe1} to dynamic random environments. 
A number of one-dimensional examples are given. 
In some cases, the sign of the speed can be determined.

\vspace{0.5cm}
\noindent
{\it Acknowledgment.} The authors are grateful to Luca Avena, Frank Redig and 
Florian V\"ollering for fruitful discussions.

\vspace{0.2cm}
\end{abstract}

\newpage

\section{Introduction}
\label{sec:intro}


\subsection{Background}
\label{subsec:RWDRE}

\emph{Random walk in random environment} (RWRE) has been an active area of research for 
more than three decades. Informally, RWRE's are random walks in discrete or continuous
space-time whose transition kernels or transition rates are not fixed but are random 
themselves, constituting a random environment. Typically, the law of the random environment 
is taken to be translation invariant. Once a realization of the random environment is 
fixed, we say that the law of the random walk is \emph{quenched}. Under the quenched law, 
the random walk is Markovian but not translation invariant. It is also interesting to 
consider the quenched law averaged over the law of the random environment, which is called 
the \emph{annealed law}. Under the annealed law, the random walk is not Markovian but 
translation invariant. For an overview on RWRE, we refer the reader to Zeitouni~\cite{Ze1,Ze2},
Sznitman~\cite{Sz1,Sz2}, and references therein.

In the past decade, several models have been considered in which the random environment 
itself evolves in time. These are referred to as \emph{random walk in dynamic random 
environment} (RWDRE). By viewing time as an additional spatial dimension, RWDRE can be 
seen as a special case of RWRE, and as such it inherits the difficulties present in RWRE 
in dimensions two or higher. However, RWDRE can be harder than RWRE because it is an interpolation 
between RWRE and homogeneous random walk, which arise as limits when the dynamics is slow, 
respectively, fast. For a list of mathematical papers dealing with RWDRE, we refer the 
reader to Avena, den Hollander and Redig~\cite{AvdHoRe2}. Most of the literature on RWDRE 
is restricted to situations in which the space-time correlations of the random environment 
are either absent or rapidly decaying. 

One paper in which a milder space-time mixing property is considered is Avena, den Hollander 
and Redig~\cite{AvdHoRe1}, where a law of large numbers (LLN) is derived for a class of 
one-dimensional RWDRE's in which the role of the random environment is taken by an 
\emph{interacting particle system} (IPS) with configuration space
\be{}
\Omega:=\{0,1\}^{\Z}.
\ee 
\vspace{-0.3cm}
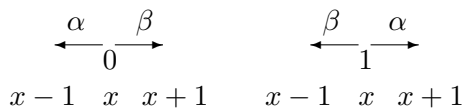
\begin{figure}[htbp]
\begin{center}
\setlength{\unitlength}{0.4cm}
\begin{picture}(12,2)(-6,1)
\put(-7.8,0){$x-1$}
\put(-4.7,0){$x$}
\put(-3.4,0){$x+1$} 
\put(-4.7,1.2){$0$}
\put(-4.3,2){\vector(1,0){1.6}}
\put(-4.7,2){\vector(-1,0){1.6}}
\put(-3.6,2.5){$\beta$}
\put(-5.9,2.5){$\alpha$}

\put(0.7,0){$x-1$}
\put(3.8,0){$x$} 
\put(5.1,0){$x+1$} 
\put(3.8,1.2){$1$}
\put(4.2,2){\vector(1,0){1.6}}
\put(3.8,2){\vector(-1,0){1.6}}
\put(4.8,2.5){$\alpha$}
\put(2.6,2.5){$\beta$}

\end{picture}
\end{center}
\vspace{0.2cm}
\caption{Jump rates of the $(\alpha,\beta)$-walk on top of a hole ($=0$), respectively, 
a particle ($=1$).}
\label{fig-jumprates}
\end{figure}

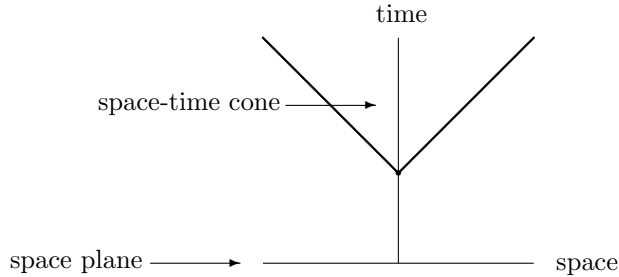
\begin{figure}[htbp]
\begin{center}
\setlength{\unitlength}{0.3cm}
\begin{picture}(12,10)(-8,1)
\put(-6,0){\line(12,0){12}} 
\put(0,0){\line(0,10){10}}
{\thicklines 
\qbezier(0,4)(3,7)(6,10) 
\qbezier(0,4)(-3,7)(-6,10)
}
\put(-11,0){\vector(1,0){4}}
\put(-5,7){\vector(1,0){4}}
\put(0,4){\circle*{.25}}
\put(-13.3,6.8){\small\mbox{space-time cone}}
\put(-1,10.7){\small\mbox{time}}
\put(7,-.3){\small\mbox{space}}
\put(-17.2,-.2){\small\mbox{space plane}}
\end{picture}
\end{center}
\vspace{0.5cm}
\caption{Cone-mixing property: asymptotic independence of states inside a space-time 
cone from states inside a space plane.}
\label{fig-cone}
\end{figure}

\noindent
In their paper, the random walk starts at $0$ and has transition rates as in Fig.~\ref{fig-jumprates}: 
on a \emph{hole} (i.e., on a $0$) the random walk has rate $\alpha$ to jump one unit to the left and rate $\beta$ 
to jump one unit to the right, while on a \emph{particle} (i.e., on a $1$) the rates are 
reversed (w.l.o.g.\ it may be assumed that $0<\beta<\alpha<\infty$, so that the random 
walk has a drift to the left on holes and a drift to the right on particles). 
Hereafter, we will refer to this model as the $(\alpha,\beta)$-model. 
The LLN is proved under the assumption that the IPS satisfies a space-time mixing property called 
\emph{cone-mixing} (see Fig.~\ref{fig-cone}), which means that the states inside a space-time cone are almost 
independent of the states in a space plane far below this cone. The proof uses a 
regeneration scheme originally developed by Comets and Zeitouni~\cite{CoZe} for RWRE 
and adapted to deal with RWDRE. This proof can be easily extended to $\Z^d$, $d \geq 2$, 
with the appropriate corresponding notion of cone-mixing.

\vspace{0.5cm}


\subsection{Elliptic vs.\ non-elliptic}
\label{subsec:ellnonell}

The original motivation for the present paper was to study the $(\alpha,\beta)$-model in
the limit as $\alpha\to\infty$ and $\beta\downarrow 0$. In this limit, which we will refer 
to as the $(\infty,0)$-model, the walk is almost a deterministic functional of the IPS; 
in particular, it is non-elliptic. The challenge was to find a way to deal with the \emph{lack 
of ellipticity}. As we will see in Section~\ref{sec:model}, our set-up will be rather general 
and will include the $(\alpha,\beta)$-model, the $(\infty,0)$-model, as well as various 
other models. Examples of papers that deal with non-elliptic (actually, deterministic) 
RW(D)RE's are Madras~\cite{Ma1} and Matic~\cite{Mat1}, where a recurrence vs.\ transience 
criterion, respectively, a large deviation principle are derived.

In the RW(D)RE literature, ellipticity assumptions play an important role. In the static case, 
RWRE in $\Z^d$, $d \geq 1$, is called \emph{elliptic} when, almost surely w.r.t.\ the 
random environment, all the rates are \emph{finite} and there is a basis $\{e_i\}_{1 \leq i
\leq d}$ of $\Z^d$ such that the rate to go from $x$ to $x+e_i$ is \emph{positive}
for $1 \leq i \leq d$. It is called \emph{uniformly elliptic} when these rates are \emph{bounded away from infinity}, 
respectively, \emph{bounded away from zero}. In \cite{CoZe}, in order to take advantage 
of the mixing property assumed on the random environment, it is important to have uniform ellipticity not necessarily in all directions, but in at least one direction in which the random walk is transient. One way to state this ``uniform directional ellipticity'' in a way that encompasses also the dynamic setting is to require the existence of a deterministic time $T>0$ and a vector $e \in \Z^d$ such that the quenched probability for the random walk to
displace itself along $e$ during time $T$ is uniformly positive for almost every realization of the random environment.
This is satisfied by the $(\alpha,\beta)$-model for $e=0$ and any $T>0$. This model is also transient 
(indeed, non-nestling) in the time direction, which enables the use of the cone-mixing property of \cite{AvdHoRe1}.
In the case of the $(\infty,0)$-model, however, there are in general no such $T$ and $e$. 
For example, when the random environment is a spin-flip system with bounded flip rates, 
any fixed space-time position has positive probability of being unreachable by the random walk. 
For all such models, the approach in \cite{AvdHoRe1} fails.

In the present paper, in order to deal with the possible lack of ellipticity we require a 
different space-time mixing property for the dynamic random environment, which we call 
\emph{conditional cone-mixing}. Moreover, as in \cite{CoZe} and \cite{AvdHoRe1}, we must 
require the random walk to have a tendency to stay inside space-time cones. Under these 
assumptions, we are able to set up a regeneration scheme and prove a LLN. Our result 
includes the LLN for the $(\alpha,\beta)$-model in \cite{AvdHoRe1}, the $(\infty,0)$-model 
for at least two subclasses of IPS's that we will exhibit, as well as models that are 
intermediate, in the sense that they are neither uniformly elliptic in any direction, 
nor deterministic as the $(\infty,0)$-model.


\subsection{Outline}
\label{subsec:outline}

The rest of the paper is organized as follows. In Section~\ref{sec:motiv} we discuss, 
still informally, the $(\infty,0)$-model and the regeneration strategy. This section 
serves as a motivation for the formal definition in Section~\ref{sec:model} of the 
class of models we are after, which is based on three \textit{structural assumptions}. 
Section~\ref{sec:results} contains the statement of our LLN under four \textit{hypotheses}, 
and a description of two classes of one-dimensional IPS's that satisfy these hypotheses 
for the $(\infty,0)$-model, namely, spin-flip systems with bounded flip rates 
that either are in Liggett's $M < \epsilon$ regime, or have finite range 
and a small enough ratio of maximal/minimal flip rates. 
Section \ref{sec:prep} contains preparation material, 
given in a general context, that is used in the proof of the LLN given in 
Section~\ref{sec:proofmainthm}. In Section~\ref{sec:proofthm2} we verify our hypotheses 
for the two classes of IPS's described in Section~\ref{sec:results}. We also obtain a 
criterion to determine the sign of the speed in the LLN, via a comparison with independent 
spin-flip systems. Finally, in  Section~\ref{sec:oex}, we discuss how to adapt the proofs 
in Section \ref{sec:proofthm2} to other models, namely, generalizations of the $(\alpha,
\beta)$-model and the $(\infty,0)$-model, and mixtures thereof. 
We also give an example where our hypotheses fail. The examples in our paper are all one-dimensional, even though our LLN is valid 
in $\Z^d$, $d\geq 1$.


\section{Motivation}
\label{sec:motiv}


\subsection{The $(\infty,0)$-model}
\label{subsec:infty0}

Let
\be{IPSdef}
\xi := (\xi_t)_{t \geq 0} \quad \mbox{ with } \quad \xi_t:= \big(\xi_t(x)\big)_{x\in\Z}
\ee 
be a c\`adl\`ag Markov process on $\Omega$. We will interpret $\xi$ by saying that at time $t$ site $x$ contains either 
a \emph{hole} ($\xi_t(x)=0$) or a \emph{particle} ($\xi_t(x)=1$). 
Typical examples are interacting particle systems on $\Omega$, such as independent spin-flips and simple exclusion.

Suppose that we run the $(\alpha,\beta)$-model on $\xi$ with 
$0 < \beta \ll 1 \ll \alpha < \infty$. Then the behavior of the random walk is as follows. 
Suppose that $\xi_0(0)=1$ and that the walk starts at $0$. The walk rapidly moves to the 
first hole on its right, typically before any of the particles it encounters manages to flip
to a hole. When it arrives at the hole, the walk starts to rapidly jump back and forth between 
the hole and the particle to the left of the hole: we say that it sits in a \emph{trap}. If
$\xi_0(0)=0$ instead, then the walk rapidly moves to the first particle on its left, where 
it starts to rapidly jump back and forth in a trap. In both cases, before moving away from 
the trap, the walk typically waits until one or both of the sites in the trap flip. If only 
one site flips, then the walk typically moves in the direction of the flip until it hits a 
next trap, etc. If both sites flip simultaneously, then the probability for the walk to sit 
at either of these sites is close to $\tfrac12$, and hence it leaves the trap in a direction 
that is close to being determined by an independent fair coin. 

The limiting dynamics when $\alpha\to\infty$ and $\beta\downarrow 0$ can be obtained from 
the above description by removing the words ``rapidly, ``typically'' and ``close to''. 
Except for the extra Bernoulli($\tfrac12$) random variables needed to decide in which 
direction to go to when both sites in a trap flip simultaneously, the walk up to time $t$ 
is a deterministic functional of $(\xi_s)_{0 \leq s \leq t}$. In particular, if 
$\xi$ changes only by single-site flips, then apart from the first jump 
the walk is completely deterministic. Since the walk spends all of its time in traps where it 
jumps back and forth between a hole and a particle, we may imagine that it lives on the 
edges of $\Z$. We implement this observation by associating with each edge its left-most 
site, i.e., we say that the walk is at $x$ when we actually mean that it is jumping back 
and forth between $x$ and $x+1$.

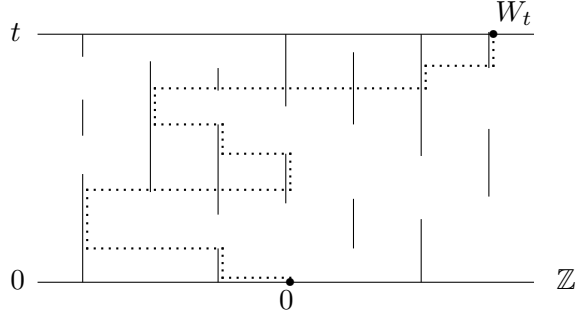
\begin{figure}[hbtp]
\vspace{0.5cm}
\begin{center}
\setlength{\unitlength}{0.3cm}
\begin{picture}(20,10)(0,0)

\put(0,0){\line(1,0){22}} 
\put(0,11){\line(1,0){22}}

\put(2,0){\line(0,1){4.8}} 
\put(2,6.5){\line(0,1){1.6}}
\put(2,10){\line(0,1){1}}

\put(5,4){\line(0,1){5.8}}

\put(8,0){\line(0,1){1.5}} 
\put(8,3){\line(0,1){4}}
\put(8,8.5){\line(0,1){1}}
\put(11,3.5){\line(0,1){2.2}}
\put(11,7.8){\line(0,1){3.2}}
\put(14,1.5){\line(0,1){2.2}}
\put(14,7){\line(0,1){3.2}} 
\put(17,0){\line(0,1){2.8}}
\put(17,5.6){\line(0,1){5.4}}
\put(20,3.8){\line(0,1){3}}
\put(20,9.5){\line(0,1){1.6}}

{\thicklines 
\qbezier[9](11.2,0.2)(9.7,0.2)(8.2,0.2)
\qbezier[4](8.2,0.2)(8.2,0.85)(8.2,1.5)
\qbezier[18](8.2,1.5)(5.2,1.5)(2.2,1.5)
\qbezier[7](2.2,1.5)(2.2,2.8)(2.2,4.1)
\qbezier[27](2.2,4.1)(6.7,4.1)(11.2,4.1)
\qbezier[5](11.2,4.1)(11.2,4.8)(11.2,5.7)
\qbezier[9](11.2,5.7)(9.7,5.7)(8.2,5.7)
\qbezier[4](8.2,5.7)(8.2,6.3)(8.2,7)
\qbezier[9](8.2,7)(6.7,7)(5.2,7)
\qbezier[5](5.2,7)(5.2,7.8)(5.2,8.6)
\qbezier[36](5.2,8.6)(11.2,8.6)(17.2,8.6)
\qbezier[3](17.2,8.6)(17.2,9.1)(17.2,9.6)
\qbezier[9](17.2,9.6)(18.7,9.6)(20.2,9.6)
\qbezier[4](20.2,9.6)(20.2,10.3)(20.2,11)
}


\put(10.7,-1.2){$0$}\put(-1.2,10.7){$t$}
\put(-1.2,-.3){$0$}
\put(20.3,11.6){$W_t$}
\put(11.2,0){\circle*{.35}}
\put(20.2,11){\circle*{.35}}
\put(23,-0.3){$\mathbb{Z}$}

\end{picture}
\end{center}
\caption{\small The vertical lines represent the presence of particles. The dotted line is the path of the $(\infty,0)$-walk.} \label{Inftyzero}
\end{figure}

Let
\be{walkdef}
W:=(W_t)_{t\geq 0}
\ee 
denote the random walk path. By the description above,  $W$ is c\`adl\`ag and
\be{eq:form0}
W_t \text{ is a function of } \big((\xi_s)_{0 \leq s \leq t} , Y\big),
\ee
where $Y$ is a sequence of i.i.d.\ Bernoulli($\tfrac12$) random variables independent of $\xi$. Note that $W$ also 
has the following three properties:
\begin{itemize}
\item[(1)] 
For any fixed time $s$, the increment $W_{s+t}-W_s$ is found by applying the same function in \eqref{eq:form0} to the environment
shifted in space and time by $(W_s,s)$ and an independent copy of $Y$; in particular, the pair $(W_t, \xi_t)$ is Markovian.
\item[(2)] 
Given that $W$ stays inside a space-time cone until time $t$, $(W_s)_{0\leq s \leq t}$ is a 
functional only of $Y$ and of the states in $\xi$ up to time $t$ inside a slightly larger cone, obtained by by adding all neighboring sites to the right.
\item[(3)]
Each jump of the path follows the same mechanism as the first jump, i.e.,
$W_t - W_{t-}$ is computed using the same rules as those for $W_0$ but applied to the environment shifted in space and time by $(W_{t-},t)$.
\end{itemize}

\noindent
The reason for emphasizing these properties will become clearer in Section~\ref{subsec:regen}.


\subsection{Regeneration}
\label{subsec:regen}

The cone-mixing property that is assumed in \cite{AvdHoRe1} to prove the LLN for the 
$(\alpha,\beta)$-model can be loosely described as the requirement that all the states 
of the IPS inside a space-time cone opening upwards depend weakly on the states inside
a space plane far below the tip (recall Fig.~\ref{fig-cone}). Let us give a rough idea 
of how this property can lead to \emph{regeneration}. Consider the event that the walk 
stands still for a long time. Since the jump times of the walk are independent of the 
IPS, so is this event. During this pause, the environment around the walk is allowed 
to mix, which by the cone-mixing property means that by the end of the pause all the 
states inside a cone with a tip at the space-time position of the walk are almost 
independent of the past of the walk. If thereafter the walk stays confined to the cone, 
then its future increments will be almost independent of its past, and so we get an 
approximate regeneration. Since in the $(\alpha,\beta)$-model there is a uniformly 
positive probability for the walk to stay inside a space-time cone with a large 
enough inclination, we see that this regeneration strategy can indeed be made to work. 

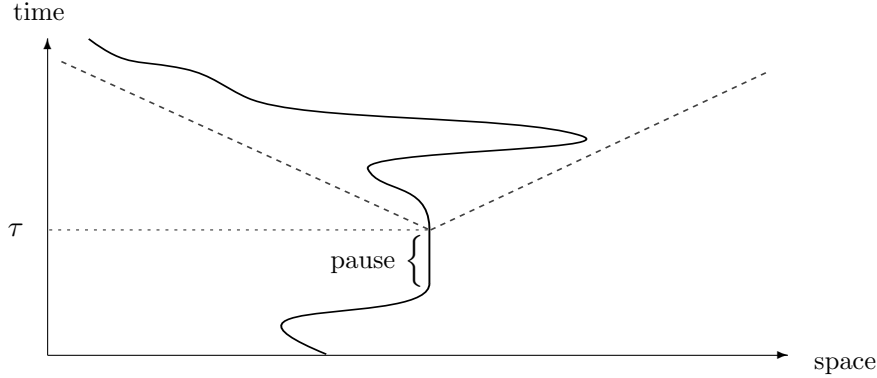
\begin{figure}
   \begin{picture}(1,0.6)
    \put(70,0){\includegraphics[height= 120pt,width=280pt]{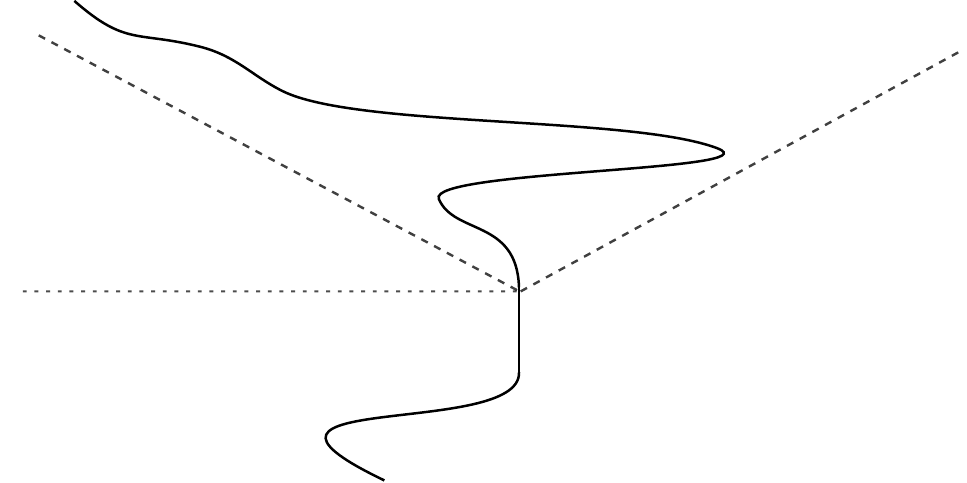}}
  \end{picture}
\put(75,0){\vector(1,0){280}}
\put(75,0){\vector(0,1){120}}
\put(62,127){\small\mbox{time}}
\put(365,-5){\small\mbox{space}}
\put(182,33){\small\mbox{pause}}
\put(210,33){\mbox{$\Big\{$}}
\put(60,45){$\tau$}
\caption{\small Regeneration at time $\tau$.} \label{fig-reg}
\end{figure}

For the actual proof of the LLN in \cite{AvdHoRe1}, cone-mixing must be more carefully 
defined. For technical reasons, there must be some uniformity in the decay of correlations 
between events in the space-time cone and in the space plane. This uniformity holds, for
instance, for any spin-flip system in the $M<\epsilon$ regime (Liggett~\cite{Li}, Section 
I.3), but not for the exclusion process or the supercritical contact process. Therefore the 
approach outlined above works for the first IPS, but not for the other two. 

There are three properties of the $(\alpha,\beta)$-model that make the above heuristics 
plausible. First, to be able to apply the cone-mixing property relative to the space-time 
position of the walk, it is important that the pair (IPS,walk) is Markovian and that the 
law of the environment as seen from the walk at any time is comparable to the initial law. 
Second, there is a uniformly positive probability for the walk to stand still for a long 
time and afterwards stay inside a space-time cone. Third, once the walk stays inside a space-time 
cone, its increments depend on the IPS only through the states inside that cone. Let us 
compare these observations with what happens in the $(\infty,0)$-model. Property (1) from 
Section~\ref{subsec:infty0} gives us the Markov property, while property (2) gives us the 
measurability inside cones. As we will see, when the environment is translation-invariant, 
property (3) implies absolute continuity of the law of the environment as seen from the walk at any positive time with respect to its 
counterpart at time zero. Therefore, as long as we can make sure that the walk has a tendency 
to stay inside space-time cones (which is reasonable when we are looking for a LLN), the 
main difference is that the event of standing still for a long time is not independent of 
the environment, but rather is a \emph{deterministic} functional of the environment. 
Consequently, it is not at all clear whether cone-mixing is enough to allow for regeneration. 
On the other hand, the event of standing still is local, since it only depends on the states 
of the two neighboring sites of the trap where the walk is pausing. For many IPS's, the 
observation of a local event will not affect the weak dependence between states that are 
far away in space-time. Hence, if such IPS's are cone-mixing, then states inside a space-time 
cone remain almost independent of the initial configuration even when we condition on seeing 
a trap for a long time.

Thus, under suitable assumptions, the event ``standing still for a long time'' is a 
candidate to induce regeneration. In the $(\alpha,\beta)$-model this event does not 
depend on the environment whereas in the $(\infty,0)$-model it is a deterministic 
functional of the environment. If we put the $(\alpha,\beta)$-model in the form (\ref{eq:form0}) 
by taking for $Y$ two independent Poisson processes with rates $\alpha$ and $\beta$, then 
we can restate the previous sentence by saying that in the $(\alpha,\beta)$-model the 
regeneration-inducing event depends only on $Y$, while in the $(\infty,0)$-model it depends 
only on $\xi$. We may therefore imagine that, also for other models of the type 
\eqref{eq:form0} and that share properties (1)--(3), it will be possible to find more 
general regeneration-inducing events that depend on both $\xi$ and $Y$ in a non-trivial 
manner. This motivates our setup in Section~\ref{sec:model}.


\section{Model setting}
\label{sec:model}

So far we have mostly been discussing RWDRE driven by an IPS. However, there are convenient 
constructions of IPS's on richer state spaces (such as graphical representations) 
that can facilitate the construction of the regeneration-inducing events mentioned in 
Section~\ref{subsec:regen}. We will therefore allow for more general Markov processes 
to represent the dynamic random environment $\xi$. Notation is set up in 
Section~\ref{subsec:notation}.  Section~\ref{subsec:mainass} contains the three 
structural assumptions that define the class of models we will consider.


\subsection{Notation and setup}
\label{subsec:notation}

Let $\N = \{1,2,\ldots \}$ be the set of natural numbers, and $\N_0 := \N \cup \{0\}$. 
Let $E$ be a Polish space and $\xi:=(\xi_t)_{t \geq 0}$ a Markov process 
with state space $E^{\Z^d}$ where $d \in \N$. Let $Y:=(Y_n)_{n \in \N}$ be an i.i.d.\ 
sequence of random elements independent of $\xi$. For $I \subset [0,\infty)$, 
abbreviate $\xi_{I}:= (\xi_u)_{u \in I}$, and analogously for $Y$. The joint law 
of $\xi$ and $Y$ when $\xi_0 = \eta \in E^{\Z^d}$ will be denoted by 
$\mathbb{P}_{\eta}$. For $n \in \N$, put $\mathscr{Y}_n:=\sigma(Y_{[1,n]})$. Let $\mathscr{F}_0 := \sigma(\xi_0)$ and, 
for $t>0$, $\mathscr{F}_t := \sigma(\xi_{[0,t]}) \vee \mathscr{Y}_{\lceil t \rceil}$. 

For $t \ge 0$ and $x \in \Z^d$, let $\theta_t$ and $\theta_x$ be the 
time-shift and space-shift operators given by
\be{}
\theta_t (\xi,Y) := \big((\xi_{t+s})_{s \geq 0},
(Y_{\lfloor t\rfloor + n})_{n\in\N}\big),
\qquad 
\theta_x (\xi,Y) := \big((\theta_x \xi_t)_{t \geq 0},
(Y_n)_{n\in\N}\big),
\ee
where $\theta_x \xi_t(y) = \xi_t(x+y)$.
In the sequel, whether $\theta$ is a 
time-shift or a space-shift operator will always be clear from the index.

We assume that $\xi$ is translation-invariant, i.e., $\theta_x\xi$ has under $\mathbb{P}_\eta$ 
the same distribution as $\xi$ under $\mathbb{P}_{\theta_x\eta}$. We also assume 
the existence of a (not necessarily unique) translation-invariant equilibrium distribution 
$\mu$ for $\xi$, and write $\mathbb{P}_{\mu}(\cdot):=\int \mu(d\eta)\,\mathbb{P}_{\eta}(\cdot)$ 
to denote the joint law of $\xi$ and $Y$ when $\xi_0$ is drawn from $\mu$.
 
The random walk will be denoted by $W=(W_t)_{t\ge0}$, and we will write $\bar\xi:=(\bar\xi_t)_{t\geq 0}$ 
to denote the \emph{environment process as seen from $W$}, i.e., $\bar\xi_t := \theta_{W_t}\xi_t$. 
Let $\bar\mu_t$ denote the law  of $\bar \xi_t$ under $\mathbb{P}_{\mu}$. 
We abbreviate $\bar\mu := \bar\mu_0$. Note that $\bar \mu = \mu$ when $\mathbb{P}_{\mu}(W_0 = 0)=1$.

For $m>0$ and $R\in\N_0$, define the $R$-enlarged $m$-cone by
\be{}
\begin{aligned}
C_R(m) &:= \big\{(x,t)\in\Z^d\times [0,\infty)\colon\,\|x\| \leq mt + R\big\},
\end{aligned}
\ee
where $\| \cdot \|$ is the $L^1$ norm. Let $\mathscr{C}_{R,t}(m)$ be 
the $\sigma$-algebras generated by the states of $\xi$ up to time $t$ inside $C_R(m)$. 


\subsection{Structural assumptions}
\label{subsec:mainass}

We will assume that $W$ is random translation of a random walk starting at $0$.
More precisely, we assume that $Z = (Z_t)_{t \ge 0}$ is a c\`adl\`ag $\mathscr{F}$-adapted $\Z^d$-valued process 
with $Z_0=0$ $\mathbb{P}_{\bar\mu}$-a.s. such that
\begin{equation}
\label{assumpWZ}
W_t = W_0 + \theta_{W_0}Z_t \quad \forall \; t \ge 0.
\end{equation}
We also assume that $W_0 \in \Z^d$ and depends on $\xi$ and $Y$ only through $\xi_0$, i.e., 
\be{assumpX0}
\mathbb{P}_{\mu}(W_0=x \mid \mathscr{F}_{\infty}) = \mathbb{P}_{\mu}(W_0=x \mid \xi_0) \; \text{ a.s. } 
\forall \; x \in \Z^d.
\ee
Under these assumptions, $(W_t-W_0)_{t\ge0}$ has under $\mathbb{P}_{\mu}$ the same distribution as $Z$ under $\mathbb{P}_{\bar \mu}$.
In what follows we make \emph{three structural assumptions} on $Z$:
\begin{itemize}
\item[(A1)] {\bf (Additivity)}\\
For all $n \in \N$,
\begin{equation}
\label{eqA1}
\begin{array}{rcl}
(Z_{t+n}-Z_n)_{t \ge 0} = \theta_{Z_n}\theta_n Z \quad \mathbb{P}_{\bar \mu} \text{-a.s.}
\end{array}
\end{equation}

\item[(A2)] {\bf (Locality)}\\
For $m>0$, let $\mathcal{D}_m := \{\|Z_t\| \le mt \,\,\forall\,
t\geq 0\}$. Then there exists $R \in \N_0$ such that,  $\forall$ $m >0$, both $\mathcal{D}_m$ and $(1_{\mathcal{D}_m}Z_t)_{t \ge 0}$ are measurable w.r.t.\ 
$\mathscr{C}_{R,\infty}(m) \vee \mathscr{Y}_{\infty}$.
\item[(A3)]{\bf (Homogeneity of jumps)}\\
For all $n\in\N$ and $x \in \Z^d$,
\begin{equation}
\mathbb{P}_{\bar \mu}\left(Z_n-Z_{n-} = x \mid \xi_{[0,n]}, 
Z_{[0,n)}\right) = \mathbb{P}_{\theta_{Z_{n-}}\xi_n}\big(W_0 = x\big)
\qquad \mathbb{P}_{\bar \mu} \text{-a.s.}
\end{equation}

\end{itemize}
These properties are analogues of properties (1)--(3) of the $(\infty,0)$-model mentioned in 
Section~\ref{subsec:infty0}, with the difference that we only require them to hold at integer times;
this will be enough as our proof relies on integer-valued regeneration times.
We also assume the `extra randomness' $Y$ to be split independently among time intervals of length $1$; 
for example, in the case of the $(\infty,0)$-model,
each $Y_n$ would \emph{not} be a Bernoulli($\frac12$) random variable but a whole \emph{sequence} of such variables instead. 
This is discussed in detail in Section \ref{sec:formaldef}.

\noindent
Another remark: assumption (A3) might seem strange since many random walk models have no deterministic jumps, 
which is indeed the case for the examples described in Section \ref{sec:results}. 
Note however that, in this case, (A3) severely restricts $W_0$, implying $W_0=0$ a.s.\ when $\xi$ is started from $\theta_{Z_{n-}}\xi_n$.
Furthermore, our main theorem (Theorem \ref{mainthm} below) is not restricted to this situation and includes also cases with deterministic jumps. For example,
one could modify the $(\infty,0)$-walk to jump exactly at integer times. Additional examples with deterministic jumps are described in item 4 of Section~\ref{sec:oex}.
The relevance of assumption (A3) is in showing that the law of the environment as seen by the RW after any jump is 
absolutely continuous w.r.t.\ the law after the first jump; this is done in Lemma \ref{lem1} below.


\section{Main results}
\label{sec:results}

Theorems~\ref{mainthm} and \ref{thm:cond} below are the main results of our paper. 
Theorem~\ref{mainthm} in Section~\ref{subsec:LLN} is our LLN. Theorem \ref{thm:cond} 
in Section~\ref{subsec:examp} verifies the hypotheses in this LLN for the $(\infty,0)$-model
in two classes of one-dimensional IPS's. For these classes some more information is 
available, namely, convergence in $L^p$, $p\geq 1$, and a criterion to determine the 
sign of the speed.


\subsection{Law of large numbers}
\label{subsec:LLN}

In order to develop a regeneration scheme for a random walk subject to assumptions 
(A1)--(A3) based on the heuristics discussed in Section~\ref{subsec:regen}, we need
 suitable regeneration-inducing events. In the \emph{four hypotheses} stated below,
  these events appear as a sequence $(\Gamma_L)_{L\in \N}$ such that, 
for a certain fixed $m \in (0,\infty)$ and $R$ as in (A2), 
$\Gamma_L\in\mathscr{C}_{R,L}(m)\vee\mathscr{Y}_L$ for all $L \in \N$.
\begin{itemize}
\item[(H1)] {\bf (Determinacy)}\\
On $\Gamma_L$, $Z_t = 0$ for all $t \in [0,L]$ $\mathbb{P}_{\bar\mu}$-a.s. 
\item[(H2)] {\bf (Non-degeneracy)}\\
For $L$ large enough, there exists a $\gamma_L>0$ such that $\mathbb{P}_{\eta}(\Gamma_L)
\geq\gamma_L$ for $\bar\mu$-a.e.\ $\eta$.
\item[(H3)] {\bf (Cone constraints)}\\
Let $\mathcal{S}:=\inf\{t > 0\colon\,\|Z_t\| > mt\}$. 
Then there exist $a \in (1, \infty)$, $\kappa_L \in (0,1]$ and $\psi_L \in [0,\infty)$ 
such that, for $L$ large enough and $\bar \mu$-a.e.\ $\eta$,
\be{}
\begin{array}{lll}
&\text{(1)} &\mathbb{P}_{\eta}(\theta_L \mathcal{S} = \infty \mid \Gamma_L ) 
\geq \kappa_L,\\[0.2cm]
&\text{(2)} &\mathbb{E}_{\eta} \left[1_{\{\theta_L \mathcal{S} < \infty\}} 
\left(\theta_L\mathcal{S}\right)^a \mid \Gamma_L \right] \leq \psi_L^a .
\end{array}
\ee
\item[(H4)] {\bf (Conditional cone-mixing)}\\
There exists a sequence of non-negative numbers  $(\Phi_L)_{L \in \N}$ 
satisfying $\lim_{L\to\infty}\kappa_L^{-1}\Phi_L=0$ such that, for $L$ large enough and for $\bar\mu$-a.e. $\eta$, \be{condconemix}
\left \lvert \mathbb{E}_{\eta} \left(\theta_L f  \mid \Gamma_L \right) 
- \mathbb{E}_{\bar \mu}(\theta_L f \mid \Gamma_L )\right\rvert 
\leq \Phi_L\,\| f\|_{\infty} \qquad \forall\,f\in\mathscr{C}_{R,\infty}(m), f \geq 0.
\ee
\end{itemize}

We are now ready to state our LLN.

\bt{mainthm}
Under assumptions {\rm (A1)--(A3)} and hypotheses {\rm (H1)--(H4)}, there exists a 
$w\in\mathbb{R}^d$ such that
\be{}
\lim_{t\to\infty} t^{-1}\,W_t = w \qquad \mathbb{P}_{\mu}-a.s.
\ee
\et

\noindent
\textbf{Remark 1:} Hypothesis (H4) above without the conditioning on $\Gamma_L$ in 
\eqref{condconemix} and with constant $\kappa_L$ is the same as the cone-mixing condition used in Avena, den Hollander 
and Redig \cite{AvdHoRe1}. There, $W_0=0$ $\mathbb{P}_{\mu}$-a.s., so that $\bar \mu = \mu$.

\noindent
\textbf{Remark 2:} Theorem \ref{mainthm} provides no information about the value of $w$, not 
even its sign when $d=1$. Understanding the dependence of $w$ on model parameters is in general 
a highly non-trivial problem.


\subsection{Examples}
\label{subsec:examp}

We next describe two classes of one-dimensional IPS's for which the $(\infty,0)$-model satisfies 
hypotheses (H1)--(H4). Further details will be given in Section~\ref{sec:proofthm2}. In both 
classes, $\xi$ is a spin-flip system in $\Omega=\{0,1\}^{\Z}$ with bounded and translation-invariant 
single-site flip rates. We may assume that the flip rates at the origin are of the form
\be{eq1}
c(\eta) = \left\{\begin{array}{ll}
c_0 + \lambda_0 p_0(\eta) & \text{ if } \; \eta(0)=1,\\
c_1 + \lambda_1 p_1(\eta) & \text{ if } \; \eta(0)=0,\\
\end{array}\right.
\qquad \eta \in \Omega,
\ee
for some $c_i, \lambda_i \geq 0$ and $p_i\colon\,\Omega\to [0,1]$, $i=0,1$.

\medskip\noindent
{\bf Example 1:} $c(\cdot)$ is in the $M<\epsilon$ regime (see Liggett~\cite{Li}, Section I.3). 

\medskip\noindent
{\bf Example 2:} $p(\cdot)$ has finite range and $(\lambda_0 + \lambda_1)/(c_0+c_1)<\lambda_c$, 
where $\lambda_c$ is the critical infection rate of the one-dimensional contact process with the same range.

\bt{thm:cond}
Consider the $(\infty,0)$-model. Suppose that $\xi$ is a spin-flip system with flip rates 
given by {\rm (\ref{eq1})}. Then for Examples $1$ and $2$ there exist a version of $\xi$ 
and events $\Gamma_L\in \mathscr{C}_{R,L}(m)\vee\mathscr{Y}_L$, $L\in\N$, 
satisfying hypotheses {\rm (H1)--(H4)}. Furthermore, the convergence in Theorem {\rm \ref{mainthm}} 
holds also in $L^p$ for all $p \geq 1$, and
\be{signofw}
\begin{array}{ll}
w \ge \frac{c_0 + \lambda_0}{c_1 + c_0 + \lambda_0}(c_1-c_0 - \lambda_0) 
& \text{ if } c_1 \ge c_0 + \lambda_0, \\[0.2cm]
w \le - \frac{c_1 + \lambda_1}{c_0 + c_1 + \lambda_1}(c_0-c_1 - \lambda_1) 
& \text{ if } c_0 \ge c_1 + \lambda_1. \\
\end{array}
\ee
\et

For independent spin-flip systems (i.e., when $\lambda_0 = \lambda_1 = 0$), \eqref{signofw} shows
 that $w$ is positive, zero or negative when the density $c_1/(c_0+c_1)$ is, 
respectively, larger than, equal to or smaller than $\tfrac12$. The criterion 
for other $\xi$ is obtained by comparison with independent spin-flip systems.

We expect hypotheses (H1)--(H4) to hold for a very large class of IPS's and walks. For each 
choice of IPS and walk, the \emph{verification} of hypotheses (H1)--(H4) constitutes a 
\emph{separate problem}. Typically, (H1)--(H2) are immediate, (H3) requires some work, 
while (H4) is hard.

Additional models will be discussed in Section~\ref{sec:oex}. We will consider generalizations 
of the $(\alpha,\beta)$-model and the $(\infty,0)$-model, namely, \emph{internal noise} models and 
\emph{pattern} models, as well as mixtures of them. 
The verification of (H1)--(H4) will be analogous to the two examples discussed above and 
will not be carried out in detail.

This concludes the motivation and the statement of our main results. The remainder of the 
paper will be devoted to the proofs of Theorems~\ref{mainthm} and \ref{thm:cond}, with the 
exception of Section~\ref{sec:oex}, which contains additional examples and remarks.


\section{Preparation}
\label{sec:prep}

The aim of this section is to prove two propositions (Propositions~\ref{prop:avrt} and 
\ref{prop:adfun} below) that will be needed in Section~\ref{sec:proofmainthm} to prove 
the LLN. In Section~\ref{sec:alln} we deal with approximate laws 
of large numbers for general discrete- or continuous-time random walks in $\mathbb{R}^d$. 
In Section~\ref{sec:adfun} we specialize to additive functionals of a Markov chain 
whose transition kernel satisfies a certain absolute-continuity property.


\subsection{Approximate law of large numbers}
\label{sec:alln}

This section contains two fundamental facts that are the basis of our proof of the LLN. 
They deal with the notion of an approximate law of large numbers. 

\begin{definition}\label{def:av}
Let $W=(W_t)_{t \geq 0}$ be a random process in $\mathbb{R}^d$ with $t\in\N_0$ or 
$t\in [0,\infty)$. For $\varepsilon \geq 0$ and $v \in \mathbb{R}^d$, we say that $W$ 
has an $\varepsilon$-approximate asymptotic velocity $v$, written $W \in AV(\varepsilon,v)$, 
if
\begin{equation}
\limsup_{t \to \infty} \left \| \frac{W_t}{t} - v \right \| \leq \varepsilon \quad 
\text{ a.s.}
\end{equation}
\end{definition}
\noindent
We take $\|\cdot\|$ to be the $L_1$-norm. A simple 
observation is that $W$ a.s.\ has an asymptotic velocity if and only if for every 
$\varepsilon > 0$ there exists a $v_{\varepsilon} \in \mathbb{R}^d$ such that $W 
\in AV(\varepsilon, v_{\varepsilon})$. In this case $\lim_{\varepsilon \downarrow 0}
v_{\varepsilon}$ exists and is equal to the asymptotic velocity.


\subsubsection{First key proposition: skeleton approximate velocity}
\label{sec:skelav}

The following proposition gives conditions under which an approximate velocity for the process 
observed along a random sequence of times implies an approximate velocity for the full process.

\begin{proposition}
\label{prop:avrt} 
Let $W$ be as in Definition~{\rm \ref{def:av}}. Set $\tau_0 := 0$, let $(\tau_k)_{k \in \N}$ be 
an increasing sequence of random times in $(0,\infty)$ (or $\N$) with 
$\lim_{k \to \infty}\tau_k = \infty$ a.s.\ and put $X_k := (W_{\tau_k},\tau_k) \in 
\mathbb{R}^{d+1}$, $k \in \N_0$. Suppose that the following hold:\\
(i) There exists an $m > 0$ such that 
\begin{equation}
\limsup_{k \to \infty} \sup_{s \in (\tau_k,\tau_{k+1}]} 
\left \| \frac{ W_s - W_{\tau_k}}{s-\tau_k}\right \| \leq m \; \text{ a.s.}
\end{equation}
(ii) There exist $v \in \mathbb{R}^d$, $u>0$ and $\varepsilon \ge 0$ such that $X \in 
AV(\varepsilon, (v,u))$.\\
Then $W \in AV((3m+1)\varepsilon/u, v/u)$.
\end{proposition}

\bpr
First, let us check that \emph{(i)} implies 
\begin{equation}
\label{sublin}
\limsup_{t \to \infty}\frac{\|W_t\|}{t} \le m \;\; \text{ a.s.}
\end{equation}
Suppose that
\begin{equation}
\label{ineq2}
\limsup_{k \to \infty} \sup_{s > \tau_k} 
\left\|\frac{W_s-W_{\tau_k}}{s-\tau_k}\right\| \leq m \quad \text{a.s.}
\end{equation}
Since, for every $k$ and $t > \tau_k$,
\begin{equation}
\left\|\frac{W_t}{t}\right\| \leq \frac{\|W_{\tau_k}\|}{t}
+\left\|\frac{W_t-W_{\tau_k}}{t-\tau_k}\right\|
\left|1-\frac{\tau_k}{t}\right| \leq \frac{\|W_{\tau_k}\|}{t}
+\sup_{s>\tau_k}\left\|\frac{W_s-W_{\tau_k}}{s-\tau_k}\right\|
\left|1-\frac{\tau_k}{t}\right|,
\end{equation}
(\ref{sublin}) follows from (\ref{ineq2}) by letting $t\to\infty$ followed by $k\to\infty$.

To check (\ref{ineq2}), define, for $k \in \N_0$ and $l \in \N$,
\begin{equation}
m(k,l) :=  \sup_{s \in (\tau_k, \tau_{k+l}]}
\left\|\frac{W_s - W_{\tau_k}}{s-\tau_k}\right\| \; \text{ and } \;
m(k,\infty) :=  \sup_{s > \tau_k}\left\|\frac{W_s - W_{\tau_k}}{s-\tau_k}\right\| = \lim_{l \to \infty}m(k,l).
\end{equation}
Using the fact that $(x_1+x_2)/(y_1+y_2)\leq (x_1/y_1)\vee(x_2/y_2)$ for all $x_1, x_2 
\in \mathbb{R}$ and $y_1, y_2 > 0$, we can prove by induction that
\begin{equation}
\label{ineq3}
m(k,l) \leq \max\{m(k,1), \dots, m(k+l-1,1)\}, \qquad l\in\N.
\end{equation}
Fix $\varepsilon>0$. By \emph{(i)}, a.s.\ there exists a $k_{\varepsilon}$ such that $m(k,1) \leq 
m+\varepsilon$ for $k>k_{\varepsilon}$. By (\ref{ineq3}), the same is true for $m(k,l)$ for 
all $l \in\N$, and therefore also for $m(k,\infty)$. 
Since $\varepsilon$ is arbitrary, (\ref{ineq2}) follows.

Let us now proceed with the proof of the proposition. Assumption $(ii)$ implies that, a.s.,
\begin{equation}
\label{eq:rt1}
\limsup_{k \to \infty}\left\| \frac{W_{\tau_k}}{k} - v \right\| 
\leq \varepsilon \;\;\; \text{ and } \;\;\; 
\limsup_{k \to \infty}\left| \frac{\tau_k}{k} - u \right| \leq \varepsilon.
\end{equation}
For $t \ge 0$, let $k_t$ be the (random) non-negative 
integer such that
\begin{equation}
\label{eq:rt2}
\tau_{k_t} \le t < \tau_{k_t+1}.
\end{equation}
Since $\tau_1 < \infty$ a.s., $k_t>0$ for large enough $t$. From (\ref{eq:rt1}) and (\ref{eq:rt2}) we deduce that
\begin{equation}
\label{eq:rt3}
\limsup_{t \to \infty}\left| \frac{t}{k_t} - u \right| \leq 
\varepsilon \quad \text{ and so } \quad 
\limsup_{t \to \infty}\left| \frac{t}{k_t} - \frac{\tau_{k_t}}{k_t} \right| \le 2\varepsilon.
\end{equation}
For $t$ large enough we may write
\begin{align}
\label{ineq:avrt1}
\left\| \frac{uW_t}{t}-v \right\| \leq 
& \frac{\|W_t\|}{t}\left|u - \frac{t}{k_t} \right| 
+ \left\|\frac{W_t - W_{\tau_{k_t}}}{k_t}\right\| 
+ \left\| \frac{W_{\tau_{k_t}}}{k_t}-v \right\|\nonumber \\
\leq & \frac{\|W_t\|}{t}\left|u - \frac{t}{k_t} \right| 
+ \sup_{s \in (\tau_{k_t},\tau_{k_t+1}]}\left\|\frac
{W_s - W_{\tau_{k_t}}}{s-\tau_{k_t}} \right\| 
\left| \frac{t - \tau_{k_t}}{k_t}\right| + \left\| \frac{W_{\tau_{k_t}}}{k_t}-v \right\|,
\end{align}
from which we obtain the conclusion by taking the limsup as $t\to\infty$ in (\ref{ineq:avrt1}), 
using (i), (\ref{sublin}), (\ref{eq:rt1}) and (\ref{eq:rt3}), and then dividing by $u$.
\epr


\subsubsection{Conditions for the skeleton to have an approximate velocity}
\label{sec:condskelav}

The following lemma states sufficient conditions for a discrete-time process to have 
an approximate velocity. It will be used in the proof of Proposition \ref{prop:adfun} below.

\bl{lemma:aplln}
Let $X=(X_k)_{k\in\N_ 0}$ be a sequence of random vectors in $\mathbb{R}^d$ with joint
law $P$ such that $P(X_0 =0)=1$. Suppose that there exist a probability measure $Q$ on $\mathbb{R}^d$ and numbers 
$\phi \in [0,1)$, $a>1$, $K>0$ with $\int_{\mathbb{R}^d} \|x\|^a\,Q(dx) \leq K^a$, such 
that, $P$-a.s. for all $k\in\N_0$,\\
(i) $\left \lvert P(X_{k+1}-X_k\in A \mid X_0,\dots, X_k) - Q(A) \right \rvert \leq \phi$ 
for all $A$ measurable;\\
(ii) $E[\|X_{k+1}-X_k\|^a|X_0,\dots,X_k] \leq K^a$.\\
Then
\be{eq:genaplln}
\limsup_{n \to \infty} \left\lVert \frac{X_n}{n} - v \right \rVert 
\leq 2K \phi^{(a-1)/a} \qquad P\text{-a.s.,}
\ee
where $v=\int_{\mathbb{R}^d} x\,Q(dx)$. In other words, $X \in AV(2K\phi^{(a-1)/a},v)$.
\el

\bpr 
The proof is an adaptation of the proof of Lemma 3.13 in \cite{CoZe}; we include it here 
for completeness. With regular conditional probabilities, we can, using $(i)$, couple $P$ 
and $Q^{\otimes\N_0}$ according to a standard splitting representation (see e.g.\ 
Berbee~\cite{Be}). More precisely, on an enlarged probability space we can construct random 
variables 
\be{}
(\Delta_k,V_k,R_k)_{k\in\N}
\ee 
such that
\begin{itemize}
\item[(1)] 
$(\Delta_k)_{k\in\N}$ is an i.i.d.\ sequence of Bernoulli($\phi$) random variables.
\item[(2)] 
$(V_k)_{k\in\N}$ is an i.i.d.\ sequence of random vectors with law $Q$.
\item[(3)] 
$(\Delta_l)_{l \geq k}$ is independent of $(\Delta_l,V_l,R_l)_{0\leq l<k}$, $R_k$.
\item[(4)]
Setting $\hat{X}_0:=0$ and, for $k \in \N_0$, $\hat{X}_{k+1}-\hat{X}_k:= (1-\Delta_k)V_k +\Delta_k R_k$, then $\hat{X}$ is equal in distribution to $X$.
\item[(5)] 
Setting $\mathcal{G}_k := \sigma(\Delta_l,V_l,R_l\colon\,0\leq l \leq k)$, then $E[ f(V_k) \mid \mathcal{G}_{k-1}]$ is measurable w.r.t.\ 
$\sigma(\hat{X}_l\colon\,0\leq l\leq k-1)$ for any Borel nonnegative function $f$.
\end{itemize}

Using (4), we may write
\be{eq5}
\frac{X_n}{n} \; \,{\buildrel d \over =}\, \; \frac{\hat{X}_n}{n} = \frac{1}{n} \sum_{k=1}^n V_k 
- \frac{1}{n} \sum_{k=1}^n \Delta_k V_k 
+ \frac{1}{n}\sum_{k=1}^n \Delta_k R_k.
\ee
As $n\to\infty$, the first term on the r.h.s.\ converges a.s.\ to $v$ 
by the LLN for i.i.d.\ random variables. By H\"older's inequality, the norm of the second term 
is at most
\be{eq2}
\left(\frac{1}{n} \sum_{k=1}^n \Delta_k \right)^{(a-1)/a} 
\left(\frac{1}{n} \sum_{k=1}^n \|V_k\|^a\right)^{1/a}, 
\ee
which, by (1) and (2), converges a.s.\ as $n\to\infty$ to
\be{}
\phi^{(a-1)/a} \left(\int_{\mathbb{R}^d} \|x\|^a\,Q(dx)\right)^{1/a}
\leq K \phi^{(a-1)/a}.
\ee

To control the third term, put $R^*_k := E[R_k\mid\mathcal{G}_{k-1}]$. Since $\|\Delta_k R_k\| \leq \|\hat{X}_{k+1}-\hat{X}_k\|$, 
using (1), (3), (4), (5) and (ii), we get
\be{eq3}
\phi E[\|R_k\|^a \mid \mathcal{G}_{k-1}] 
=  E[\Delta_k \|R_k\|^a \mid \mathcal{G}_{k-1}] 
\leq E[\|\hat{X}_{k+1}-\hat{X}_k\|^{a}\mid\mathcal{G}_{k-1}] \leq K^a.
\ee
Combining \eqref{eq3} with Jensen's inequality, we obtain
\be{eq4}
\|R^*_k\| \leq E\big[\|R_k\|^{a} \mid \mathcal{G}_{k-1}\big]^{1/a} 
\leq \frac{K}{\phi^{1/a}},
\ee
so that
\be{eq6}
\left\|\frac{1}{n} \sum_{k=1}^n \Delta_k R^*_k \right\| 
\leq \frac{K}{\phi^{1/a}} 
\left(\frac{1}{n} \sum_{k=1}^n \Delta_k \right) \xrightarrow[n \to \infty]{} K \phi^{(a-1)/a}.
\ee

Now fix $y \in \mathbb{R}^d$ and put
\be{}
M^y_n := \sum_{k=1}^n \frac{\Delta_k}{k} \langle R_k-R^*_k, y \rangle.
\ee
where $\langle \cdot, \cdot \rangle$ denotes the usual inner product. Then $(M^y_n)_{n\in\N_0}$ is 
a $(\mathcal{G}_n)_{n\in\N_0}$-martingale whose quadratic variation is
\be{}
\langle M^y \rangle_n = \sum_{k=1}^{n} \frac{\Delta_k}{k^2}\langle R_k - R^*_k,y \rangle^2.
\ee
By the Burkholder-Gundy inequality and (\ref{eq3}--\ref{eq4}), we have
\be{}
\begin{aligned}
E\left[\sup_{n\in\N}|M^y_n|^{a \wedge 2}\right] 
&\leq C\,E \left[\langle M^y \rangle_{\infty}^{(a \wedge 2)/2}\right]\\ 
&\leq C\,E \left[\sum_{k=1}^{\infty} \frac{\Delta_k}{k^{a \wedge 2}}
\left| \langle R_k-R^*_k,y \rangle \right |^{a \wedge 2}\right] \leq C\,\|y\|^{a \wedge 2} 
K^{a \wedge 2},
\end{aligned}
\ee
where $C$ is a positive constant that may change after each inequality. This implies that $M^y_n$ is uniformly integrable
 and therefore converges a.s.\ as $n\to\infty$. Kronecker's lemma then gives
\be{}
\lim_{n\to\infty} \frac{1}{n} \sum_{k=1}^n \Delta_k \langle R_k-R^*_k, y \rangle 
= 0 \, \text{ a.s.}
\ee
Since $y$ is arbitrary, this in turn implies that
\be{eq7}
\lim_{n\to\infty} \frac{1}{n} \sum_{k=1}^n \Delta_k(R_k-R^*_k) = 0 \; \text{ a.s.}
\ee
Therefore, by \eqref{eq6} and \eqref{eq7}, the limsup of the norm of the last term in the r.h.s.\ of \eqref{eq5} is also bounded by $K\phi^{(a-1)/a}$, which finishes the proof.
\epr


\subsection{Additive functionals of a discrete-time Markov chain}
\label{sec:adfun}


\subsubsection{Notation}

Let $\mathcal{X} = (\mathcal{X}_n)_{n \in \N_0}$ be a time-homogeneous Markov chain in the canonical space equipped 
with the time-shift operators $(\theta_n)_{n \in \N_0}$. For $n \ge 1$, put $\mathcal{F}_n := 
\sigma(\mathcal{X}_{[1,n]})$ (note that $\mathcal{X}_0 \notin \mathcal{F}_{\infty}$) and let $P_\chi$ denote the law of $(\mathcal{X}_n)_{n 
\in \N_0}$ when $\mathcal{X}_0 = \chi$. Fix an initial measure $\nu$ and suppose that, for any nonnegative $f \in \mathcal{F}_{\infty}$, 
\begin{equation}
\label{sup:abscont}
P_\nu(E_{\mathcal{X}_n}[f] \in \cdot) \ll P_\nu(E_{\mathcal{X}_0}[f] \in \cdot),
\end{equation}
where $P_\nu := \int \nu(d\chi) P_\chi$.

Let $Z = (Z_n)_{n \in \N_0}$ be a $\Z^d$-valued $\mathcal{F}$-adapted process that is an \emph{additive functional} of $(\mathcal{X}_n)_{n \in \N}$, i.e., 
$Z_0=0$ and, for any $k \in \N_0$, 
\begin{equation}\label{defadditivefunc}
(Z_{k+n}-Z_{k})_{n \in \N_0} = \theta_k Z \qquad P_\nu \text{-a.s.}
\end{equation}

We are interested in finding random times $(\tau_k)_{k \in \N_0}$ such 
that $(Z_{\tau_k},\tau_k)_{k \in \N_0}$ satisfies the hypotheses of 
Lemma~\ref{lemma:aplln}. In the Markovian setting it makes sense to look for 
$\tau_k$ of the form
\begin{equation}
\label{eq:recrt}
\tau_0 = 0, \qquad \tau_{k+1} = \tau_k + \theta_{\tau_k}\tau, \quad k\in\N_0,
\end{equation}
where $\tau$ is a random time.

Condition (i) of Lemma~\ref{lemma:aplln} is a ``decoupling condition''. It states 
that the law of an increment of the process depends weakly on the previous increments. Such a 
condition can be enforced by the occurrence of a ``decoupling event'' under which 
the increments of $(Z_{\tau_k},\tau_k)_{k \in \N_0}$ lose dependence. In this 
setting, $\tau$ is a time at which the decoupling event is observed.


\subsubsection{Second key proposition: approximate regeneration times}

Proposition~\ref{prop:adfun} below is a consequence of Lemma~\ref{lemma:aplln} and 
is the main result of this section. It will be used together with Proposition~\ref{prop:avrt} 
to prove the LLN in Section~\ref{sec:proofmainthm}. It gives a way to construct $\tau$ 
when the decoupling event can be detected by ``probing the future'' with a stopping time.

For a random variable $\mathcal{T}$ taking values in $\N_0 \cup \{\infty\}$, we define the 
\emph{image} of $\mathcal{T}$ by $\mathcal{I}_{\mathcal{T}}:= \{n \in \N\colon\, 
P_{\nu}(\mathcal{T} = n)>0\}$, and its closure under addition by $\bar{\mathcal{I}}_{\mathcal{T}}
:= \{n \in \N\colon\,\exists\,\,l \in \N,\,i_1,\ldots,i_l \in 
\mathcal{I}_{\mathcal{T}}\colon\, n = i_1 + \dots + i_l\}$.
Note that $\mathcal{I}_{\mathcal{T}} = \emptyset$ if and only if $\mathcal{T} \in \{0,\infty\}$ a.s..

\begin{proposition}
\label{prop:adfun}
Let $\mathcal{T}$ be a stopping time for the filtration $\mathcal{F}$
taking values in $\N\cup\{\infty\}$. Put $D := \{\mathcal{T} = \infty\}$ and suppose that
the following properties hold:\\
(i) For every $n \in \bar{\mathcal{I}}_{\mathcal{T}}$ there exists a $D_n \in \mathcal{F}_n$ 
such that 
\begin{equation*}
D \cap \theta_n D = D_n \cap \theta_n D \quad \; P_{\nu} \text{-a.s.}
\end{equation*}
(ii) There exist numbers $\rho \in (0,1]$, $a>1$, $C>0$, $m>0$ and $\phi \in [0,1)$ such that, $P_{\nu}$-a.s.,
\begin{equation*}
\begin{array}{rl}
\text{(a)} & P_{\mathcal{X}_0}\left(D \right) \geq \rho, \\
\text{(b)} & E_{\mathcal{X}_0}\left[1_{\{\mathcal{T} < \infty\}} \mathcal{T}^a \right] \leq C^a, \\
\text{(c)} & \text{On } D \text{, } \|Z_t\| \leq mt \text{ for all } t \in \N_0,\\
\text{(d)} & \left| E_{\mathcal{X}_0}\left[f(Z, (\theta_n\mathcal{T})_{n\in \bar{\mathcal{I}}_{\mathcal{T}}}) \mid D\right] - E_{\nu}\left[f(Z, (\theta_n\mathcal{T})_{n\in \bar{\mathcal{I}}_{\mathcal{T}}}) \mid D\right] \right| \le \phi \|f\|_{\infty} \; \forall f \ge 0 \text{ measurable.}
\end{array}
\end{equation*}
Then there exists a random time $\tau \in \mathcal{F}_{\infty}$ taking values in $\N$ 
such that, setting $\tau_k$ as in {\rm (\ref{eq:recrt})} and $X_k:= (Z_{\tau_k}, \tau_k)$, 
then $X \in AV(\varepsilon, (v,u))$ where $(v,u)=E_{\nu}\left[(Z_{\tau},\tau) \mid 
D\right]$, $u>0$ and $\varepsilon = 12(m+1)u\phi^{(a-1)/a}$. 
\end{proposition}


\subsubsection{Two further propositions}

In order to prove Proposition~\ref{prop:adfun}, we will need two further propositions (Propositions \ref{prop:adfun1} and \ref{prop:stoptimes} below). 

\begin{proposition}
\label{prop:adfun1}
Let $\tau$ be a random time measurable w.r.t.\ $\mathcal{F}_{\infty}$ taking values in 
$\N$. Put $\tau_k$ as in {\rm (\ref{eq:recrt})} and $X_k := (Z_{\tau_k},\tau_k)$. 
Suppose that there exists an event $D \in \mathcal{F}_{\infty}$ such that the following 
hold $P_{\nu}$-a.s.:\\
(i) For $n \in \mathcal{I}_{\tau}$, there exist events $H_n$ and $D_n \in \mathcal{F}_n$ 
such that
\begin{equation}
\begin{array}{rl} 
\text{(a)} & \{\tau = n\} = H_n \cap \theta_n D, \\ 
\text{(b)} & D \cap \theta_n D = D_n \cap \theta_n D.
\end{array}
\end{equation}
(ii) There exist $\phi \in [0,1)$, $K>0$ and $a>1$ such that that, on $\{P_{\mathcal{X}_0}(D)>0\}$,
\begin{equation}
\begin{array}{rl}
\text{(a)} & E_{\mathcal{X}_0}[\|X_1\|^a | D] \le K^a,\\
\text{(b)} &\left | P_{\mathcal{X}_0}\left(X_1 \in A | D \right) 
-  P_{\nu}\left(X_1 \in A | D \right) \right | \leq \phi 
\quad \forall\,A \text{ measurable.}\\
\end{array} 
\end{equation}
Then $X \in AV\big(\varepsilon, (v,u)\big)$, where $\varepsilon= 2K\phi^{(a-1)/a}$ and $(v,u)
:= E_{\nu}[X_1|D]$.
\end{proposition}

\bpr
Since $\tau < \infty$, by (i)(a) and (\ref{sup:abscont}) we must have $P_{\nu}(D)>0$.
Let $\mathcal{F}_{\tau_k}$ be the $\sigma$-algebra of the events $B \in \mathcal{F}_{\infty}$ 
such that, for all $n \in \N$, there exists 
$B_n \in \mathcal{F}_n$ with $B \cap \{\tau_k = n\} = B_n \cap \{\tau_k = n\}$. We will 
show that, $P_{\nu}$-a.s., for all $k\in\N$,
\begin{equation}
\label{eq:af1}
E_{\nu}\left[ \| \theta_{\tau_k}X_1 \|^a | \mathcal{F}_{\tau_k}\right] \le K^a
\end{equation}
and
\begin{equation}
\label{eq:af2}
\left|P_{\nu}\left(\theta_{\tau_k}X_1 \in A | \mathcal{F}_{\tau_k}\right) 
- P_{\nu}(X_1 \in A | D) \right| \leq \phi \quad \forall\,A \text{ measurable.}
\end{equation}
Then, setting $Q(\cdot):= P_{\nu}(X_1 \in \cdot | D)$ and noting that $\theta_{\tau_k}X_1 = X_{k+1} - X_k$ 
and $X_j \in \mathcal{F}_{\tau_k}$ for all $0\leq j \leq k$, we will be able to conclude since
(\ref{eq:af1}--\ref{eq:af2}) and (ii)(a) imply that the conditions of Lemma~\ref{lemma:aplln} are all satisfied.

To prove (\ref{eq:af1}--\ref{eq:af2}), first note that, using (i), one can verify by induction
that (i)(a) holds also for $\tau_k$, i.e., for every $n \in \mathcal{I}_{\tau_k}$ there 
exists $H_{k,n} \in \mathcal{F}_n$ such that
\begin{equation}
\{\tau_k = n \} =  H_{k,n} \cap \theta_n D \qquad P_\nu \text{-a.s.}
\end{equation}
Take $B \in \mathcal{F}_{\tau_k}$ and a measurable nonnegative function $f$, and write
\begin{align}
\label{eq:af3}
E_{\nu}\left[ 1_B \theta_{\tau_k}f(X_1)\right] 
& = \sum_{n \in \mathcal{I}_{\tau_k}}E_{\nu}
\left[ 1_{B \cap \{\tau_k = n\}}\theta_n f(X_1)\right] 
= \sum_{n \in \mathcal{I}_{\tau_k}} 
E_{\nu}\left[1_{B_n \cap H_{k,n}} 
\theta_n\big(1_D f(X_1)\big) \right] \nonumber \\ 
& = \sum_{n \in \mathcal{I}_{\tau_k}} E_{\nu}\left[1_{B_n \cap H_{k,n}} 
P_{\mathcal{X}_n}(D)E_{\mathcal{X}_n}\left[f(X_1) | D \right]\right].
\end{align}
Noting that $P_{\nu}(B) = \sum_{n \in \mathcal{I}_{\tau_k}} E_{\nu}\left[1_{B_n \cap H_{k,n}} P_{\mathcal{X}_n}(D) \right]$, 
obtain (\ref{eq:af1}) by taking $f(x) = \|x\|^a$ and using (ii)(a) together 
with (\ref{sup:abscont}). For (\ref{eq:af2}), choose $f = 1_A$, subtract $P_{\nu}(B)
E_{\nu}\left[ f(X_1) | D \right]$ from (\ref{eq:af3}) and use (ii)(b).
\epr

\noindent

\begin{proposition}
\label{prop:stoptimes}
Let $\mathcal{T}$ be a stopping time as in Proposition \ref{prop:adfun} 
and suppose that conditions (ii)(a) and (ii)(b) of that proposition are satisfied. 
Define a sequence of stopping times $(T_k)_{k \in \N_0}$ as follows. 
Put $T_0 = 0$ and, for $k \in \N_0$,
\begin{equation}
\label{recursionT_k}
T_{k+1} := \left\{
\begin{array}{ll}
\infty & \;\;\;\;\;\; \text{ if } T_k = \infty \\
 T_k + \theta_{T_k}\mathcal{T} & \;\;\;\;\;\; \text{ otherwise.}
\end{array}\right.
\end{equation}
Put
\begin{equation}
N := \inf\{ k \in \N_0\colon\, T_k < \infty \text{ and } T_{k+1} = \infty\}.
\end{equation}
Then $N < \infty$ a.s.\ and there exists a constant $\varkappa = \varkappa(a,\rho) \in (0,\infty)$ such that, 
$P_\nu$-a.s.,
\begin{equation}
\label{conclem}
E_{\mathcal{X}_0}\left[T_N^a \right] \le (\varkappa C)^a.
\end{equation}
Furthermore, $\mathcal{I}_{T_N} \subset \bar{\mathcal{I}}_{\mathcal{T}}$.
\end{proposition}

\bpr
First, let us check that
\begin{equation}
\label{nfin}
P_{\mathcal{X}_0}(N \geq n) \le (1-\rho)^n.
\end{equation}
Indeed, $N \geq n$ if and only if $T_n < \infty$, so that, for $k \in \N_0$,
\begin{equation}
\label{nfinn}
P_{\mathcal{X}_0}(T_{k+1} < \infty) 
= E_{\mathcal{X}_0}\left[ 1_{\{T_k < \infty\}} 
P_{\mathcal{X}_{T_k}}(\mathcal{T} < \infty)\right] \le (1-\rho)P_{\mathcal{X}_0}(T_k < \infty),
\end{equation}
where we use (ii)(a) and the fact that (\ref{sup:abscont}) holds also with a stopping time in place of $n$. 
Clearly, (\ref{nfin}) follows from (\ref{nfinn}) by induction. In particular, $N < \infty$ a.s.

From \eqref{recursionT_k} we see that, for $0\leq k \leq n$,
\begin{equation}\label{recursion2}
T_n = T_k + \theta_{T_k}T_{n-k} \; \text{ on } \; \{T_k < \infty \}.
\end{equation}
Using (ii)(a) and (b), with the help of (\ref{sup:abscont}) again, we can a.s.\ estimate, 
for $0\leq k < n$,
\begin{align}
\label{estim}
E_{\mathcal{X}_0}\left[1_{\{T_n < \infty \}} \left|T_{k+1}-T_k \right|^a\right] 
& = E_{\mathcal{X}_0}\left[1_{\{T_{k+1} < \infty\}} \left|T_{k+1}-T_k \right|^a 
P_{\mathcal{X}_{T_{k+1}}}(T_{n-k-1}<\infty)\right] \nonumber \\ 
& \le (1-\rho)^{n-k-1} E_{\mathcal{X}_0}
\left[1_{\{T_k < \infty, \theta_{T_k}\mathcal{T} < \infty\}} 
\theta_{T_k}\mathcal{T}^a\right] \nonumber \\ 
& = (1-\rho)^{n-k-1} E_{\mathcal{X}_0}\left[1_{\{T_k < \infty\}} 
E_{\mathcal{X}_{T_k}}\left[1_{\{\mathcal{T} < \infty\}} 
\mathcal{T}^a \right] \right] \nonumber \\ 
& \le (1-\rho)^{n-k-1} C^a P_{\mathcal{X}_0}(T_k < \infty) \nonumber \\ 
& \le (1-\rho)^{n-1} C^a.
\end{align}
Now write
\begin{equation}
T_N = \sum_{k=0}^{N-1} T_{k+1}-T_k.
\end{equation}
By Jensen's inequality,
\begin{equation}
T_N^a \le N^{a-1}\sum_{k=0}^{N-1} \left| T_{k+1}-T_k \right|^a
\end{equation}
so that, by \eqref{estim},
\begin{align}
E_{\mathcal{X}_0}\left[ T_N^a \right] \le \sum_{n=1}^{\infty} n^{a-1}
\sum_{k=0}^{n-1}E_{\mathcal{X}_0}\left[ 1_{\{N=n\}}
\left|T_{k+1}-T_k \right|^a \right] \le C^a \sum_{n=1}^{\infty} n^a (1-\rho)^{n-1} \; \text{ a.s.}
\end{align}
and (\ref{conclem}) follows by taking $\varkappa = \left( \sum_{n=1}^{\infty} n^a 
(1-\rho)^{n-1} \right)^{1/a}$.

As for the claim that $\mathcal{I}_{T_N} \subset \bar{\mathcal{I}}_{\mathcal{T}}$,
write, for $n \in \N$,
\begin{equation}
\{T_N = n\} = \sum_{k=1}^{\infty} \{T_k=n, N=k\}
\end{equation}
to see that $\mathcal{I}_{T_N} \subset \bigcup_{k=1}^{\infty} \mathcal{I}_{T_k}$. Using 
(\ref{recursionT_k}), we can verify by induction that, for each $k \in \N$, $\mathcal{I}_{T_k} 
\subset \{n \in \N \colon\,\exists\,i_1,\ldots,i_k \in \mathcal{I}_{\mathcal{T}}
\colon\,n = i_1 + \dots + i_k \} \subset \bar{\mathcal{I}}_{\mathcal{T}}$, and the claim 
follows.
\epr


\subsubsection{Proof of Proposition \ref{prop:adfun}}

We can now combine Propositions \ref{prop:adfun1} and \ref{prop:stoptimes} to prove 
Proposition~\ref{prop:adfun}.

\bpr 
In the following we will refer to the hypotheses of Proposition~\ref{prop:adfun1} 
with the prefix P. For example, P(i)(a) denotes hypothesis (i)(a) in that proposition.
The hypotheses in Proposition~\ref{prop:adfun} will be referred to without a prefix. Since 
the hypotheses of Proposition~\ref{prop:stoptimes} are a subset of those of 
Proposition~\ref{prop:adfun}, the conclusions of the former are valid.

We will show that, if $\tau := t_0 + \theta_{t_0}T_N$ for a suitable $t_0 \in \N$, then $\tau$ 
satisfies the hypotheses of Proposition~\ref{prop:adfun1} for a suitable $K$. There are 
two cases. If $\mathcal{I}_{\mathcal{T}} = \emptyset$, then $T_N \equiv 0$. Choosing 
$t_0 = 1$, we basically fall in the context of Lemma~\ref{lemma:aplln}. P(i)(a) and 
P(i)(b) are trivial, (ii)(c) implies that P(ii)(a) holds with $K = (m+1)$, while P(ii)(b) 
follows immediately from (ii)(d). Therefore, we may suppose that $\mathcal{I}_{\mathcal{T}} \neq \emptyset$ 
and put $\iota := \min \mathcal{I}_{\mathcal{T}} \in \N$. 
Let $\hat C := 1 \vee (\varkappa C)$ and $t_0 := \iota \lceil \hat C \rho^{-1/a} \rceil$. 
We will show that $\tau$ satisfies the hypotheses of Proposition~\ref{prop:adfun1} 
with $K = 6\iota(m+1)\hat C\rho^{-1/a}$.

\medskip\noindent
\underline{P(i)(a)}:
First we show that this property is true for $T_N$. Indeed,
\begin{align}
\{T_N = n\} & = \sum_{k \in \N_0}\{N=k, T_k = n\} 
= \sum_{k \in \N_0}\{T_k=n, \theta_n \mathcal{T} = \infty\} \\
& = \theta_n D \cap \left( \bigcup_{k \in \N_0} \{T_k=n\} \right),
\end{align}
and $\hat H_n := \bigcup_{k \in \N_0} \{T_k=n\} \in \mathcal{F}_n$ since the $T_k$'s 
are all stopping times. Now we observe that $\{ \tau = n\} = \theta_{t_0}\{T_N = n-t_0\}$, 
so we can take $H_n := \emptyset$ if $n<t_0$ and $H_n := \theta_{t_0} \hat H_{n-t_0}$ 
otherwise.

\medskip\noindent
\underline{P(i)(b)}:
By $(i)$, it suffices to show that $\mathcal{I}_{\tau} \subset \bar{\mathcal{I}}_{\mathcal{T}}$. 
Since $t_0 \in \bar{\mathcal{I}}_{\mathcal{T}}$ (as an integer multiple of $\iota$), this 
follows from the definition of $\tau$ and the last conclusion of Proposition \ref{prop:stoptimes}.

\medskip\noindent
\underline{P(ii)(a)}:
By (ii)(c), $\|X_1\|^a  = \left(\|Z_{\tau}\|+\tau\right)^a \le \left((m+1)\tau\right)^a$ on 
$D$. Therefore, we just need to show that
\begin{equation}
\label{mombound}
E_{\mathcal{X}_0}\left[\tau^a |D\right] \le (6 \iota \hat C)^a/\rho.
\end{equation}
Now, $\tau^a \le 2^{a-1}\left(t_0^a + \theta_{t_0}T_N^a\right)$ and, by 
Proposition~\ref{prop:stoptimes} and (\ref{sup:abscont}),
\begin{equation}
E_{\mathcal{X}_0}\left[ \theta_{t_0} T_N^a\right] 
= E_{\mathcal{X}_0}\left[ E_{\mathcal{X}_{t_0}}\left[T_N^a\right]\right] \le \hat C^a.
\end{equation}
Using (ii)(a), we obtain
\begin{equation}
E_{\mathcal{X}_0}\left[ \theta_{t_0} T_N^a\right | D] \le \hat C^a / \rho.
\end{equation}
Since $t_0 \le 2 \iota \hat C \rho^{-1/a}$ and $\iota \ge 1$, (\ref{mombound}) follows.

\medskip\noindent
\underline{P(ii)(b)}:
Let $S = (S_n)_{n\in \bar{\mathcal{I}}_{\mathcal{T}}}$ with $S_n := \theta_n \mathcal{T}$. By (ii)(d), it is enough 
to show that $X_1 = (Z_{\tau}, \tau) \in \sigma(Z,S)$ a.s.. Since $Z_{\tau} = \sum_{n=0}^{\infty}
1_{\{\tau = n\}}Z_n \in \sigma(Z, \tau)$, it suffices to show that $\tau \in \sigma(S)$ a.s..
Using the definition of the $T_k$'s, we verify by induction that each $T_k$ is a.s.\ measurable in 
$\sigma(S)$. Since $N \in \sigma((T_k)_{k \in \N_0})$, both $N$ and $T_N$ are also a.s.\
in $\sigma(S)$. Therefore, a.s. $\tau \in \sigma(\theta_{t_0}S) \subset \sigma(S)$.  
 
\medskip
With all hypotheses verified, Proposition~\ref{prop:adfun1} implies that $X \in AV(\hat 
\varepsilon, (v,u))$, where $(v,u) = E_{\nu}[X_1|D]$ and $\hat \varepsilon = 2K\phi^{(a-1)/a}$. 
To conclude, observe that $u = E_{\nu}[\tau|D] \ge t_0 \ge \iota \hat C \rho^{-1/a} > 0$, 
so that $K = 6(m+1) \iota \hat C \rho^{-1/a} \le 6(m+1)u$. Therefore, $\hat \varepsilon 
\le \varepsilon$ and the proposition follows. In the case $\mathcal{I}_{\mathcal{T}} = 
\emptyset$, we conclude similarly since $u=1$ and $K=(m+1)$.
\epr


\section{Proof of Theorem~\ref{mainthm}}
\label{sec:proofmainthm}

In this section we show how to put the model defined in Section~\ref{sec:model} in the 
context of Section~\ref{sec:prep}, and we prove the LLN using Propositions~\ref{prop:avrt} 
and \ref{prop:adfun}.


\subsection{Two further lemmas}

Before we start, we first derive two lemmas (Lemmas~\ref{lem1} and \ref{lem2} below) that 
will be needed in Section~\ref{subsec:proofmainthm}. The first lemma relates the laws of 
the environment as seen from $W_n$ and from $W_0$. The second lemma is an extension of the 
conditional cone-mixing property for functions that depend also on $Y$.

\bl{lem1}
$\bar{\mu}_n\ll \bar\mu$ for all $n \in \N$.
\el

\bpr
For $t\geq 0$, let $\bar\mu_{t-}$ denote the law of $\theta_{W_{t-}}\xi_t$ under 
$\mathbb{P}_{\mu}$. First we will show that $\bar\mu_{t-} \ll \mu$. This is a consequence 
of the fact that $\mu$ is translation-invariant equilibrium, and remains true if we replace 
$W_{t-}$ by any random variable taking values in $\Z^d$. Indeed, if $\mu(A)=0$ then 
$\mathbb{P}_{\mu}(\theta_x \xi_t \in A)=0$ for every $x \in \Z^d$, so
\be{}
\begin{aligned}
\bar{\mu}_{t-}(A) = \mathbb{P}_{\mu}(\theta_{W_{t-}}\xi_t \in A) 
= \sum_{x \in \Z^d} \mathbb{P}_{\mu}(W_{t-}=x, \theta_x\xi_t \in A) = 0.
\end{aligned}
\ee
Now take $n \in \N$ and let $g_n := \frac{d\bar\mu_{n-}}{d \mu}$. For any measurable 
$f\geq 0$,
\begin{align}
\mathbb{E}_{\mu}\left[f(\theta_{W_n}\xi_n)\right] = \mathbb{E}_{\bar \mu}\left[f(\theta_{Z_n}\xi_n)\right]
= & \sum_{x \in \Z^d} \mathbb{E}_{\bar \mu}
\left[1_{\{Z_n-Z_{n-}=x\}}f(\theta_x\theta_{Z_{n-}}\xi_n)\right] \nonumber \\
 = & \sum_{x \in \Z^d} \mathbb{E}_{\bar \mu}
\left[\mathbb{P}_{\theta_{Z_{n-}}\xi_n}(W_0=x)f(\theta_x\theta_{Z_{n-}}\xi_n)\right] \nonumber \\
 = & \sum_{x \in \Z^d} \mathbb{E}_{\mu}
\left[\mathbb{P}_{\theta_{W_{n-}}\xi_n}(W_0=x)f(\theta_x\theta_{W_{n-}}\xi_n)\right] \nonumber \\ 
 = & \sum_{x \in \Z^d} \mathbb{E}_{\mu}
\left[g_n(\xi_0)\mathbb{P}_{\xi_0}(W_0=x)f(\theta_x\xi_0)\right] \nonumber \\
 = & \sum_{x \in \Z^d} \mathbb{E}_{\mu}
\left[g_n(\xi_0)1_{\{W_0=x\}} f(\theta_x\xi_0)\right] \nonumber\\
 = & \; \mathbb{E}_{\mu} \left[ g_n(\xi_0) f(\theta_{W_0} \xi_0)\right] 
\end{align}
where, for the second equality, we use (A3). 
\epr

\bl{lem2}
For $L$ large enough and for all nonnegative $f\in\mathscr{C}_{R,\infty}(m) \vee \mathscr{Y}_{\infty}$,
\be{eq:lem2}
\left\lvert\mathbb{E}_{\eta} \left[\theta_L f \mid \Gamma_L\right] 
- \mathbb{E}_{\bar \mu}[\theta_Lf\mid\Gamma_L ]\right\rvert 
\leq \Phi_L \|f\|_{\infty} \quad \text{for } \;\bar{\mu} \text{-a.e. } \eta.
\ee
\el	
\bpr
Put $f_y(\eta) = f(\eta,y)$ and abbreviate $Y^{(L)} = (Y_k)_{k > L}$. Then $\theta_L f 
= \theta_L f_{Y^{(L)}}$. Since $\Gamma_L$ depends on $Y$ only through $(Y_k)_{k \le L}$, we have
\be{}
\mathbb{E}_{\eta}[\theta_Lf\,1_{\Gamma_L} \mid Y^{(L)}] 
= \mathbb{E}_{\eta}\big[\theta_Lf_{(\cdot)}\,1_{\Gamma_L}\big] \circ (Y^{(L)}),
\ee
and \eqref{eq:lem2} follows from (H4) applied to $f_y$.
\epr


\subsection{Proof of Theorem~\ref{mainthm}}
\label{subsec:proofmainthm}
\bpr
Extend $\xi$ and $Z$ for times $t \in [-1,0]$ by taking them constant in this interval, 
and let $Y_0$ be a copy of $Y_1$ independent of $\mathscr{F}_{\infty}$. Put
\begin{equation}\label{defmkvchain}
\begin{array}{rcl}
\mathcal{X}_0 & := & \left( \xi_{[-1,0]},Z_{[-1,0]},Y_0\right), \\
\mathcal{X}_{n+1} & := & \left(\theta_{Z_n}\xi_{[n,n+1]},(Z_{t+n}-Z_n)_{0 \le t \le 1}, Y_{n+1}\right), \; n \in \N_0.
\end{array}
\end{equation}
Then $(\mathcal{X}_n)_{n \in \N_0}$ is a time-homogeneous Markov chain; 
to avoid confusion, we will denote its time-shift operator by $\bar\theta_n$.
Note that $\mathcal{F}_n = \mathscr{F}_n$ $\forall$ $n \in \N\cup\{\infty\}$ and that, for functions $f \in \mathcal{F}_{\infty}$, 
$\bar \theta_n f = \theta_{Z_n} \theta_n f$ $\forall$ $n \in \N_0$.

Fix $L \in \N$ large enough and put
\begin{equation}\label{defcurlyTL}
\mathcal{T}_L := L + 1_{\Gamma_L} \lceil \theta_L \mathcal{S} \rceil.
\end{equation}
By (\ref{eqA1}) and since $\Gamma_L \in \mathscr{F}_L$ and $Z$ is $\mathscr{F}$-adapted, 
$\mathcal{T}_L$ is an $\mathcal{F}$-stopping time and
$(Z_n)_{n \in \N_0}$ is an additive functional of $(\mathcal{X}_n)_{n \in \N}$ as in Section \ref{sec:adfun}.

Next, we will verify (\ref{sup:abscont}) for $\mathcal{X}$ and the hypotheses of Proposition~\ref{prop:adfun} for 
$Z$ and $\mathcal{T}_L$ under $\mathbb{P}_{\bar \mu}$. These hypotheses will be referred 
to with the prefix P. The notation here is consistent in the sense that parameters in 
Section~\ref{sec:model} are named according to their role in Section~\ref{sec:prep}; the 
presence/absence of a subscript $L$ indicates whether the parameter depends on $L$ or not.

\medskip\noindent
\underline{(\ref{sup:abscont})}:
Noting that, for nonnegative $f \in \mathcal{F}_{\infty}$ and $n \in \N_0$,
\begin{equation}\label{relMCwithmodel}
E_{\mathcal{X}_n} \left[ f \right] = \mathbb{E}_{\theta_{Z_n}\xi_n}\left[ f \right] \qquad \mathbb{P}_{\bar \mu} \text{-a.s.},
\end{equation}
this follows from Lemma \ref{lem1} and (\ref{assumpWZ}--\ref{assumpX0}).

\medskip\noindent
\underline{P(i)}:
We will find $D_n$ for $n\ge L$. This is enough, since both $\mathcal{I}_{\mathcal{T}_L}$ 
and $\bar{ \mathcal{I}}_{\mathcal{T}_L}$ are subsets of $[L,\infty) \cap \N$. Using (A1) and (H1), we 
may write
\begin{equation}
\begin{array}{rcl}
D & = & \Gamma_L \cap \{\|Z_{t+L}\|\le mt \; \forall \; t \ge 0\}, \\
\bar\theta_n D & = & \bar \theta_n \Gamma_L \cap \{\|Z_{t+n+L}-Z_n\|\le mt \; \forall \; t \ge 0\}.\\
\end{array}
\end{equation}
Intersecting the two above events, we get
\begin{equation}
D\cap\bar\theta_n D = \Gamma_L \cap \{\|Z_t\| \le mt 
\; \forall \; t \in [0,n]\} \cap \bar \theta_n D,
\end{equation}
i.e., P(i) holds with $D_n :=\Gamma_L \cap \{\|Z_t\| \le mt \; \forall \; 
t \in [0,n]\} \in \mathcal{F}_n$ for $n \ge L$.

For the remaining items, note that, by \eqref{relMCwithmodel},
the distribution of $(Z,\mathcal{T}_L)$ under $P_{\mathcal{X}_0}$ is $\mathbb{P}_{\bar \mu}$-a.s. the same as under $\mathbb{P}_{\xi_0}$.

\medskip\noindent
\underline{P(ii)(a)}:
Since $\{\mathcal{T}_L = \infty\} = \{\theta_L\mathcal{S} = \infty\} \cap \Gamma_L$, we get 
from (H2) and (H3)(1) that, $\mathbb{P}_{\bar \mu}$-a.s.,
\begin{equation}
 \mathbb{P}_{\xi_0}\left(\mathcal{T}_L = \infty \right) 
= \mathbb{P}_{\xi_0}\left( \theta_L \mathcal{S} = \infty \mid \Gamma_L\right)\mathbb{P}_{\xi_0}(\Gamma_L)\ge \kappa_L \gamma_L > 0,
\end{equation}
so that we can take $\rho_L := \kappa_L \gamma_L$.

\medskip\noindent
\underline{P$(ii)(b)$}:
By the definition of $\mathcal{T}_L$, we have
\begin{align}
\mathcal{T}_L^a 1_{\{\mathcal{T}_L < \infty\}} 
& = L^a 1_{\Gamma_L^c} + \left(L+  \lceil \theta_L \mathcal{S} \rceil \right)^a 
1_{\Gamma_L \cap \{\theta_L \mathcal{S} < \infty\}}\nonumber \\
& \le L^a 1_{\Gamma_L^c} + \left(L+1+\theta_L\mathcal{S} \right)^a 
1_{\Gamma_L \cap \{\theta_L \mathcal{S} < \infty\}}\nonumber \\
& \le 2^{a-1}(L+1)^a + 2^{a-1} \left((\theta_L \mathcal{S})^a
1_{ \{\theta_L \mathcal{S}<\infty\}}\right) 1_{\Gamma_L} .
\end{align}
Therefore, by (H3)(2), we get
\begin{align}
\mathbb{E}_{\xi_0}\left[ \mathcal{T}_L^a 1_{\{\mathcal{T}_L<\infty\}}\right] 
\le 2^a((L+1)^a+(1\vee\psi_L)^a) \le \left[2(L+1+1\vee\psi_L)\right]^a \quad \mathbb{P}_{\bar \mu} \text{-a.s.},
\end{align}
so that we can take $C_L := 2(L+1+1\vee\psi_L)$.

\medskip\noindent
\underline{P(ii)(c)}: 
This follows from (H1) and the definition of $\mathcal{S}$.

\medskip\noindent
\underline{P(ii)(d)}:
First note that, for any $n \in \bar{\mathcal{I}}_{\mathcal{T}_L}$, $\bar\theta_n \mathcal{T}_L \in \sigma(Z, \bar\theta_n\Gamma_L)$. 
Since $n \ge L$, on $\{\mathcal{T}_L = \infty\} = \Gamma_L \cap \{\theta_L \mathcal{S}=\infty\}$, 
$Z$, $\bar\theta_n \Gamma_L$ and $\{\theta_L \mathcal{S} = \infty\}$ are all measurable in $\theta_L(\mathscr{C}_{R,\infty}(m) \vee \mathscr{Y}_{\infty})$;
 this follows from (A2), (H1) and the assumptions on $\Gamma_L$. 
Noting that, for any two probability measures $\nu_1$, $\nu_2$ and an event $A$, 
\begin{equation}\label{compcondTV}
\| \nu_1(\cdot \mid A) - \nu_2(\cdot \mid A) \|_{TV} \le 2 \frac{\| \nu_1 - \nu_2 \|_{TV}}{\nu_1(A) \vee \nu_2(A)}
\end{equation}
where $\| \cdot \|_{TV}$ stands for total variation distance, we see that
P(ii)(d) follows from Lemma \ref{lem2} and (H3)(1) with $\phi_L := 2 \Phi_L/\kappa_L \to 0$ as $L \to \infty$ by (H4).

\medskip
Thus, for large enough $L$, we can conclude by Proposition~\ref{prop:adfun} that there 
exists a sequence of times $(\tau_k)_{k \in \N_0}$ with  $\lim_{k\to\infty}
\tau_k = \infty$ a.s.\ such that $(Z_{\tau_k},\tau_k)_{k \in \N_0} \in 
AV(\varepsilon_L, (v_L,u_L))$, where 
\begin{equation}
\label{vLuL}
\begin{array}{rcl}
v_L & = & \mathbb{E}_{\bar\mu}[Z_{\tau_1}|D], \\
u_L & = & \mathbb{E}_{\bar\mu}[\tau_1|D] > 0, \\
\varepsilon_L & = & 12(m+1)u_L \phi_L^{(a-1)/a}.
\end{array}
\end{equation}
From (\ref{vLuL}) and P(ii)(c), Proposition~\ref{prop:avrt} implies that $Z \in 
AV(\delta_L, w_L)$, where
\begin{equation}
\begin{array}{rcl}
w_L & = & v_L/u_L, \\
\delta_L & = & (3m+1)12(m+1)\phi_L^{(a-1)/a}.
\end{array}
\end{equation}
By (H4), $\lim_{L\to\infty}\delta_L = 0$. As was observed after Definition~\ref{def:av}, 
this implies that $w := \lim_{L \to \infty} w_L$ exists and that $\lim_{t \to \infty}t^{-1}
Z_t = w$ $\mathbb{P}_{\bar\mu}$-a.s., which, by (\ref{assumpWZ}--\ref{assumpX0}), implies the same for 
$W$, $\mathbb{P}_{\mu}$-a.s.
\epr

We have at this point finished the proof of our LLN. In the following sections, we will 
look at examples that satisfy (H1)--(H4). Section~\ref{sec:proofthm2} is devoted to the 
$(\infty,0)$-model for two classes of one-dimensional spin-flip systems. In Section~\ref{sec:oex} 
we discuss three additional models where the hypotheses are satisfied, and one where 
they are not.


\section{Proof of Theorem \ref{thm:cond}}
\label{sec:proofthm2}

We begin with a proper definition of the $(\infty,0)$-model in Section~\ref{sec:formaldef}, 
where we identify $Z$ and $W_0$ of Section~\ref{subsec:mainass}. 
In Section~\ref{sec:sfsbr}, we first define suitable versions of spin-flip systems with bounded 
rates. After checking assumptions (A1)--(A3), we define events 
$\Gamma_L$ satisfying (H1) and (H2) for which we then verify (H3). We also derive uniform integrability properties of 
$t^{-1}W_t$ which are the key for convergence in $L^p$ once we have the LLN. In 
Sections~\ref{sec:mlessep} and \ref{sec:subcritdep}, we specialize to particular constructions 
in order to prove (H4), which is the hardest of the four hypotheses. Section~\ref{sec:signspeed} 
is devoted to proving a criterion for positive or negative speed.


\subsection{Definition of the model}
\label{sec:formaldef}

Assume that $\xi$ is a c\`adl\`ag process with state space $E:=\{0,1\}^{\Z}$. We will
define the walk $W$ in several steps, and a monotonicity property will follow.

\subsubsection{Identification of $Z$ and $W_0$}
\label{identZW0}

First, let $Tr^+=Tr^+(\eta)$ and $Tr^-=Tr^-(\eta)$ denote the locations of the closest 
traps to the right and to the left of the origin in the configuration $\eta \in E$, i.e.,
\begin{equation}
\label{deftraps}
\begin{array}{rcl}
Tr^+(\eta) & := & \inf \{x  \in \N_0\colon\,\eta(x) = 1, \eta(x+1)=0\},\\
Tr^-(\eta) & := & \sup \{x  \in -\N_0\colon\,\eta(x) = 1, \eta(x+1)=0\},
\end{array}
\end{equation}
with the convention that $\inf\emptyset = \infty$ and $\sup\emptyset = - \infty$. For $i,j 
\in \{0,1\}$, abbreviate $\langle i,j \rangle := \{\eta \in E \colon\, \eta(0) = i, \eta(1) = j\}$. 
Let $\bar E := \langle 1,0 \rangle$, i.e., the set of all the configurations with a trap 
at the origin.

Next, we define the functional $J$ that gives the jumps in $W$. For $b \in \{0,1\}$ and 
$\eta \in E$, let
\begin{equation}
\label{defJ}
J(\eta,b) := Tr^+\left( 1_{\langle 1,1 \rangle} 
+ b1_{\langle 0,1 \rangle} \right) + Tr^-\left( 1_{\langle 0,0 \rangle} 
+ (1-b)1_{\langle 0,1 \rangle} \right),
\end{equation}
i.e., $J$ is equal to either the left or the right trap, depending on the configuration around the 
origin. In the case of an inverted trap ($\langle 0,1 \rangle$), the 
direction of the jump is decided by the value of $b$. Observe that $J=Tr^+=Tr^-=0$ when $\eta \in 
\bar E$, independently of the value of $b$. 

Let $b_0$ be a Bernoulli($\frac12$) random variable independent of $\xi$ and set
\begin{equation}\label{defW0}
W_0 = X_0 := J(\xi_0,b_0).
\end{equation}
Now let $(b_{n,k})_{n,k \in \N}$ be a double-indexed i.i.d.\ sequence 
of Bernoulli($\frac12$) r.v.'s independent of $(\xi,b_0)$.
Put $\tau_0 := 0$ and, for $k \geq 0$,
\begin{equation}
\label{defjumps}
\begin{array}{rcll}
\tau_{k+1} & := & \left\{ \begin{array}{l} \infty \\
                         \inf\left\{t> \tau_k \colon\, \left(\xi_t(X_{k}),\xi_t(X_{k}+1)\right) \neq (1,0)\right\}
                         \end{array}\right.  & \begin{array}{l} \text{if } |X_k| = \infty,\\ \text{otherwise,} \end{array} \\
X_{k+1} & := & \left\{ \begin{array}{l} X_k \\
                        X_k + J\left(\theta_{X_k}\xi_{\tau_k},b_{\lceil \tau_{k+1} \rceil, k+1}\right)\\
                        \end{array}\right.  & \begin{array}{l} \text{if } \tau_{k+1} = \infty,\\ \text{otherwise.} \end{array}
\end{array}
\end{equation}
Since $\xi$ is c\`adl\`ag, for any $k \in \N_0$ we either have $\tau_{k} = \infty$ or $\tau_{k+1} > \tau_{k}$. We define 
$(W_t)_{t\geq 0}$ as the path that jumps $X_{k+1} - X_k$ at time $\tau_{k+1}$ and is constant 
between jumps, i.e.,
\begin{equation}
\label{defW}
W_t := \sum_{k=0}^{\infty}1_{\{\tau_k \le t<\tau_{k+1}\}} X_k.
\end{equation}
With this definition, it is clear that $W_t$ is c\`adl\`ag and, by (\ref{defW0}--\ref{defjumps}),
\begin{equation}\label{Wadditive}
W_{n+t}-W_n = \theta_{W_n} \theta_n W_t \;\; \text{on } \{W_n < \infty\} \;\; \forall \; n \in \N_0, t \ge 0.
\end{equation}
Therefore, defining $Z$ by 
\begin{equation}\label{defZ}
Z_t := 1_{\{ \xi_0 \in \bar E \}} W_t, \qquad t \ge 0,
\end{equation}
we get $W_t = W_0 + \theta_{W_0}Z_t$ on $\{W_0 < \infty\}$ 
since, in this case, $\theta_{W_0}\xi_0 \in \bar E$, and $W_0=0$ on $\bar E$. 

\subsubsection{Monotonicity}
\label{subsubsec:monotonicity}

The following monotonicity property will be helpful in checking (H3). In order to state it, 
we first endow both $E$ and $D([0,\infty),E)$ with the usual partial ordering, i.e., for $\eta_1, 
\eta_2 \in E$, $\eta_1 \leq \eta_2$ means that $\eta_1(x)\leq\eta_2(x)$ for all $x \in \Z$, while, 
for $\xi^{(1)}, \xi^{(2)} \in D([0,\infty),E)$, $\xi^{(1)} \le \xi^{(2)}$ means that $\xi^{(1)}_t 
\leq \xi^{(2)}_t$ for all $t \geq 0$.

\begin{lemma}
\label{monot}
Fix a realization of $b_0$ and $(b_{n,k})_{n,k \in \N}$. If $\xi^{(1)} \le \xi^{(2)}$, then 
\begin{equation}\label{eqmonot}
W_t\left(\xi^{(1)},b_0,(b_{n,k})_{n,k \in \N}\right) \leq W_t\left(\xi^{(2)}, b_0, (b_{n,k})_{n,k \in \N}\right)
\end{equation} 
for all $t \ge 0$.
\end{lemma}
\bpr
This is a straightforward consequence of the definition. We need only to understand what 
happens when the two walks separate and, at such moments, the second walk is always to the right of 
the first.
\epr


\subsection{Spin-flip systems with bounded flip rates}
\label{sec:sfsbr}

\subsubsection{Dynamic random environment}
\label{subsubsec:DRE}

From now on we will take $\xi$ to be a single-site spin-flip system with translation-invariant 
and bounded flip rates. We may assume that the rates at the origin are of the form
\begin{equation}
\label{rates}
c(\eta) = \left\{ \begin{array}{ccl}
c_0 + \lambda_0 p_0(\eta) & \text{ when } & \eta(0) = 1, \\
c_1 + \lambda_1 p_1(\eta) & \text{ when } & \eta(0) = 0,
\end{array}
\right.
\end{equation}
where $c_i,\lambda_i > 0$ and $p_i \in [0,1]$. We assume the existence conditions of Liggett~\cite{Li}, 
Chapter I, which in our setting amounts to the additional requirement that $c(\cdot)$ has finite triple norm.
This is automatically satisfied in the $M<\epsilon$ regime or when $c(\cdot)$ has finite range.

From (\ref{rates}), we see that the IPS is stochastically dominated by the system $\xi^+$ with rates
\begin{equation}
\label{upperisf}
c^+(\eta) = \left\{ \begin{array}{ccl}
c_0 & \text{ when } & \eta(0) = 1, \\
c_1 + \lambda_1 & \text{ when } &\eta(0) = 0,
\end{array}\right.
\end{equation}
while it stochastically dominates the system $\xi^-$ with rates
\begin{equation}
\label{lowerisf}
c^-(\eta) = \left\{ \begin{array}{ccl}
c_0 + \lambda_0& \text{ when } & \eta(0) = 1, \\
c_1 & \text{ when } & \eta(0) = 0.
\end{array}
\right.
\end{equation}
These are the rates of two independent spin-flip systems with respective densities 
$\rho^+ := (c_1+\lambda_1)/\lambda^+$ and $\rho^- := c_1/\lambda^-$ where $\lambda^+ := c_0+c_1+\lambda_1$ and $\lambda^- := c_0 + \lambda_0+c_1$.
Consequently, any equilibrium for $\xi$ is stochastically dominated by $\nu_{\rho^+}$ and dominates $\nu_{\rho^-}$, 
where $\nu_\rho$ is a Bernoulli product measure with density $\rho$.

We will take as the dynamic random environment the triple $\Xi := (\xi^-,\xi, \xi^+)$ 
starting from the same initial configuration and coupled together via the basic 
(or Vasershtein) coupling, which implements 
the stochastic ordering as an a.s.\ partial ordering. More precisely, $\Xi$ is the IPS with 
state space $E^3$ whose rates are translation invariant and at the origin are given 
schematically by (the configuration of the middle coordinate is $\eta$), 
\begin{equation}\
\label{ratescoupl}
\begin{array}{ccl}
(000) & \to & \left\{\begin{array}{cl}
      (111) & c_1, \\
      (011) & c(\eta) - c_1, \\
      (001) & c_1 + \lambda_1 - c(\eta), \\
      \end{array}\right.\\
(001) & \to & \left\{\begin{array}{cl}
      (111) & c_1,\\
      (011) & c(\eta) - c_1, \\
      (000) & c_0,\\
      \end{array}\right. \\
(011) & \to & \left\{\begin{array}{cl}
      (111) & c_1, \\
      (000) & c_0, \\
      (001) & c(\eta) - c_0, \\
      \end{array}\right.\\
(111) & \to & \left\{\begin{array}{cl}
      (000) & c_0, \\
      (001) & c(\eta) - c_0, \\
      (011) & c_0 + \lambda_0 - c(\eta). \\
      \end{array}\right. \\
\end{array}
\end{equation}

\subsubsection{Verification of (A1)--(A3)}
\label{checkassumptions}

Under our assumptions, $\lim_{k \to \infty}\tau_k = \infty$ and $X_0 < \infty$ $\mathbb{P}_\mu$-a.s., as $\xi$ has bounded flip rates per site and $\mu$ dominates and is dominated by non-trivial product measures. By induction, $X_k < \infty$ a.s.\ for every $k \in \N$ as well, since the law of $\theta_{X_{k-1}}\xi_{\tau_k}$ is absolutely continuous w.r.t.\ $\mu$, which can be verified by approximating $\tau_k$ from above by times taking values in a countable set. Therefore, $W_t$ is finite for all $t \ge 0$.

Set $Y_n := (b_{n,k})_{k \in \N}$. Then $Z$ is $\mathscr{F}$-adapted as it is independent of $b_0$.
(A1) follows by \eqref{Wadditive} and \eqref{defZ}, and (A3) follows either from the recursive construction \eqref{defjumps} 
or by noting that $Z$ has no deterministic jumps and $\theta_{Z_n} \theta_n W_0 = 0$.
To verify (A2), note that $\{J=x\}$ depends on $\eta$ only through $(\eta(y))_{ y \in \{0 \wedge x, \ldots, 0 \vee x + 1\}}$
so we may take $R=1$.

\subsubsection{Definition of $\Gamma_L$ and verification of (H1)--(H3)}
\label{subsubsec:defGammaL}

Using $\Xi$, we can define the events $\Gamma_L$ by
\begin{equation}
\Gamma_L := \big\{\xi_t^{\pm}(x)=\xi_0^{\pm}(x) \; 
\forall \; t \in [0,L] \text{, } x = 0,1\big\}.
\end{equation}
Then $\Gamma_L \in \mathscr{C}_{1,L}(m)$ for any $m > 0$. 
When $\xi_0^{\pm} \in \bar E$, $\Gamma_L$ implies that there is a trap at the origin between times $0$ and $L$; since $\bar \mu(\bar E)=1$, (H1) holds. The probability of $\Gamma_L$ is positive 
and depends on $\Xi_0$ only through the states at $0$ and $1$, 
so (H2) is also satisfied.

In order to verify (H3), we will take advantage of Lemma~\ref{monot} and the stochastic domination in 
$\Xi$ to define two auxiliary processes $H^{\pm} = (H^{\pm}_t)_{t\ge0}$ which we can control and which will bound $Z$. 
This will also allow us to deduce uniform integrability properties.

In the following we will suppose that $\xi^{{\pm}}_0 \in \bar E$. Let $G_0 = U_0 := 0$ and, for $k \geq 0$,

\begin{equation}\label{defUG+}
\begin{array}{rcl}	
U_{k+1} & := & \inf \left\{t>U_k \colon\, \xi^{+}_t(G_k+1)=1 \right\}, \\
G_{k+1} & := & G_k + Tr^{+}\Big(\theta_{G_k}\xi^{+}_{U_{k+1}}\Big) \\
\end{array}
\end{equation}
and put
\begin{equation}\label{defH+}
H^{+}_t := \sum_{k=0}^{\infty} 1_{\{U_k \le t < U_{k+1}\}}G_{k+1}.
\end{equation}
Define $H^-$ analogously, using $Tr^-$ and $\xi^-$ instead and switching $1$'s to $0$'s in \eqref{defUG+}.
Then $H^+$ ($H^-$) is the process that, observing $\xi^+$ ($\xi^-$) , waits to the left of a hole (on a particle) until it flips to a 
particle (hole), and then jumps to the right (left) to the next trap. 
Therefore, by Lemma~\ref{monot} and the definition of $Z$, $H^-_t \le Z_t \le H^+_t$ $\forall$ $t \ge 0$. 
Note that $H^+$ depends only on $(\xi^+(x))_{x \ge 1}$, and analogously for $H^-$.

In the following, we will write $\Z_{\le x} : = \Z \cap (-\infty,x]$ and analogously for $\Z_{\ge x}$.

\begin{lemma}
\label{sbisf}
Fix $\rho_* \in (0,\rho^-]$ and $\rho^* \in [\rho^+,1)$. 
There exist $m,a,\psi_*$ $\in (0,\infty)$ and $\kappa_* \in (0,1)$, depending on $\rho_*$, $\rho^*$ and $\lambda^{\pm}$, such that, for any probability measure $\bar \nu$ on $\bar E$ that stochastically dominates $\nu_{\rho_*}$ on $\Z_{\le-1}$ and is dominated by $\nu_{\rho^*}$ on $\Z_{\ge 2}$,
\begin{equation}\label{sbisfeq1}
\text{(a)} \;\; \sup_{t \ge 1}\mathbb{E}_{\bar \nu}\left[e^{a \left( t^{-1}|H^{\pm}_t| \right)} \right] \le \psi_*
\end{equation}
and, setting
\begin{equation}\label{defSpm}
\mathcal{S}^{\pm}:= \inf\{t > 0 \colon\, |H^{\pm}_t| > mt\}, \qquad \widehat{\mathcal{S}}^{\pm}:= \sup\{t > 0 \colon\, |H^{\pm}_t| > mt\},
\end{equation}
then
\begin{equation}\label{sbisfeq3}
\begin{array}{rl}
\text{(b)} & \mathbb{P}_{\bar\nu}\left(\mathcal{S}^{\pm} = \infty \right) \ge \kappa_*,\\
\text{(c)} & \mathbb{E}_{\bar\nu}\left[ e^{a \widehat{\mathcal{S}}^{\pm}} \right] \le \psi_*.
\end{array}
\end{equation}
\end{lemma}

Before proving this lemma, let us see how it leads to (H3). 
We will show that there exist $m,a,\psi \in (0,\infty)$ and $\kappa \in (0,1)$
such that, for all $L \ge 1$ and $\eta \in \bar E$,
\begin{equation}\label{verifH31}
\mathbb{P}_{\eta} \left( \theta_L\mathcal{S}=\infty \mid \Gamma_L \right) \ge \kappa 
\end{equation}
and
\begin{equation}\label{verifH32}
\mathbb{E}_{\eta} \left[ e^{a (\theta_L \mathcal{S}) } 1_{\{\theta_L\mathcal{S}<\infty \}} \mid \Gamma_L \right] \le \psi,
\end{equation}
which clearly imply (H3). 

Let us verify \eqref{verifH31}.
First note that $\theta_L \mathcal{S} \ge \theta_L (\mathcal{S}^+ \wedge \mathcal{S}^-)$, 
and that the latter is nonincreasing in $(\eta(x))_{x \ge 2}$ and nondecreasing in $(\eta(x))_{x \le -1}$.
Therefore we may assume that $\eta = \eta_{01}$ which is the configuration in $\bar E$ 
with all $0$'s on $\Z_{\le-1}$ and all $1$'s on $\Z_{\ge 2}$. In this case, 
$\xi^-_L$ is distributed as $\nu_{\rho^L_0}$ on $\Z_{\le-1}$ and $\xi^+_L$ as $\nu_{\rho^L_1}$ on $\Z_{\ge 2}$,
where $\rho_0^L = \rho^-(1-e^{-\lambda^- L})$ and $\rho_1^L = \rho^+ + e^{-\lambda^+ L}(1-\rho^+)$. 
Furthermore, on $\Gamma_L$, $\xi^{\pm}_L \in \bar E$. 

Let now $m,a,\psi_*$ and $\kappa^*$ as in Lemma~\ref{sbisf} for $\rho_*:=\rho_0^1$ and $\rho^* := \rho_1^1$, and let $\bar\nu_L$ be the distribution of $\bar \eta_L \in \bar E$ given by $\xi^-_L$ on $\Z_{\le-1}$ and $\xi^+_L$ on $\Z_{\ge 2}$.
Noting that $\bar \eta_L$ is independent of $\Gamma_L$ and that $\mathcal{S}^+$ and $\mathcal{S}^-$ are independent, 
we use the previous observations, the Markov property and Lemma~\ref{sbisf}(b) to write
\begin{align}\label{compverifH31}
\mathbb{P}_{\eta}\left(\theta_L \mathcal{S} = \infty \mid \Gamma_L\right) 
& \ge \mathbb{P}_{\eta_{01}}\left(\theta_L(\mathcal{S}^+ \wedge \mathcal{S}^-) = \infty \mid \Gamma_L\right) \nonumber \\
& = \mathbb{E}_{\eta_{01}}\left[ 1_{\Gamma_L} \mathbb{P}_{\bar \eta_L}\left(\mathcal{S}^+ \wedge \mathcal{S}^- = \infty\right)\right]\mathbb{P}_{\eta_{01}}\left(\Gamma_L\right)^{-1} \nonumber \\
& = \mathbb{P}_{\bar\nu_L}\left(\mathcal{S}^+ = \infty\right)\mathbb{P}_{\bar\nu_L}\left(\mathcal{S}^- = \infty\right) \ge \kappa_*^2 \in (0,1),
\end{align}
and we may take $\kappa := \kappa_*^2$.
For \eqref{verifH32}, note now that, when finite, $\theta_L\mathcal{S} < \theta_L(\widehat{\mathcal{S}}^+ \vee \widehat{\mathcal{S}}^-)$ 
and the latter is is nondecreasing in $(\eta(x))_{x \ge 2}$ and nonincreasing in $(\eta(x))_{x \le -1}$. 
Therefore we may again assume $\eta= \eta_{01}$ and write, using Lemma~\ref{sbisf}(c),
\begin{align}\label{compverifH32}
\mathbb{E}_{\eta}\left[\theta_L \left(e^{a \mathcal{S}}1_{\{\mathcal{S}=\infty\}}\right) \mid \Gamma_L\right] 
& \le \mathbb{E}_{\eta_{01}}\left[ \theta_L e^{a \left( \widehat{\mathcal{S}}^+ + \widehat{\mathcal{S}}^- \right) } \mid \Gamma_L\right] \nonumber \\
& = \mathbb{E}_{\bar\nu_L}\left[e^{a \widehat{\mathcal{S}}^+} \right] \mathbb{E}_{\bar\nu_L}\left[e^{a \widehat{\mathcal{S}}^-} \right] \le \psi_*^2 \in (0,\infty),
\end{align}
and we can take $\psi := \psi^2_*$. All that is left to do is to prove Lemma~\ref{sbisf}.

\bpr[Proof of Lemma~\ref{sbisf}]
By symmetry, it is enough to prove (a)--(c) for $H^+$.
Since $H^+$, $\mathcal{S}^+$ and $\widehat{\mathcal{S}}^+$ are monotone, we may assume that
$\xi^+$ has rates $\lambda^+ \rho^*$ to flip from holes to particles
and $\lambda^+ (1-\rho^*)$ from particles to holes and starts from $\nu_{\rho^*}$, which is the equilibrium measure.
In this case, the increments $G_{k+1}-G_k$ are i.i.d.\ Geom($1-\rho^*$), 
and $U_{k+1}-U_k$ are i.i.d.\ Exp($\lambda^+ \rho^*$), independent from $(G_k)_{k \in \N_0}$. 
Therefore, $H^+$ is a c\`adl\`ag L\'evy process and $H^+_1$ has an exponential moment, so (a) promptly follows.
Moreover, $H^+$ satisfies a large deviation estimate of the type
\begin{equation}\label{pr:sbisfeq1}
\mathbb{P}_{\nu_{\rho^*}}\left(\exists \; s>t \text{ such that } H^+_s > ms\right) 
\leq K_1 e^{-K_2 t} \; \text{ for all } t > 0,
\end{equation}
where $m$, $K_1$ and $K_2$ are functions of ($\rho^*$, $\lambda^+$), which proves (c). 
In particular, $\widehat{\mathcal{S}}^+ < \infty$ a.s., which implies that $\mathbb{P}_{\nu_{\rho^*}}(H^+_s \le m(s+n^*) \; \forall \; s \ge 0) \ge \frac12$ for some $n^*$ large enough;
then
\begin{align}\label{pr:sbisfeq2}
\mathbb{P}_{\nu_{\rho^*}}\left(\mathcal{S}^+ = \infty \right) 
& \ge \mathbb{P}_{\nu_{\rho^*}}\left( H^+_{n^*} =0, H^+_{n^*+s} - H^+_{n^*} \le m(s+n^*) \; \forall \; s \ge 0 \right) \nonumber\\
& = \mathbb{P}_{\nu_{\rho^*}}\left( H^+_{n^*} =0 \right)\mathbb{P}_{\nu_{\rho^*}}\left(H^+_s  \le m(s+n^*) \; \forall \; s \ge 0 \right) =: \kappa_* >0,
\end{align}
proving (b).
\epr

\subsubsection{Uniform integrability}
\label{subsubsec:ui}

The following corollary implies that, for systems given by (\ref{rates}), $(t^{-1}|W_t|^p)_{t \geq 1}$ 
is uniformly integrable for any $p \geq 1 $, so that, whenever we have a LLN, the convergence holds 
also in $L^p$.

\begin{corollary}
\label{cor:ui}
Let $\xi$ be a spin-flip system with rates as in \eqref{rates}, starting from equilibrium. Then 
$(t^{-1}W_t)_{t \geq 1}$ is bounded in $L^p$ for all $p \ge 1$.
\end{corollary}

\bpr
The claim for $Z$ under $\mathbb{P}_{\bar \mu}$ follows from Lemma~\ref{sbisf}(a) 
by noting that $\bar \mu$ stochastically dominates $\nu_{\rho^-}$ on $\Z_{\le-1}$ and is dominated by $\nu_{\rho^+}$ on $\Z_{\ge 2}$; this can be verified noting that $W_0\ge0$ corresponds to finding particles to the left of $W_0$, and $W_0\le0$ to holes to its right. The same for $W$ follows from (\ref{assumpWZ}--\ref{assumpX0}) since $W_0$ has exponential moments under $\mathbb{P}_{\mu}$.
\epr

We still need to verify (H4). This will be done in Sections~\ref{sec:mlessep} and 
\ref{sec:subcritdep} below. 
As $\kappa$ in \eqref{verifH31} could be taken independently of $L$ for (H3),
we only need $\lim_{L \to \infty}\Phi_L = 0$ in (H4).


\subsection{Example 1: $M < \epsilon$}
\label{sec:mlessep}

We recall the definition of $M$ and $\epsilon$ for a translation-invariant spin-flip system:
\begin{eqnarray}
\label{defm}
M &:=& \sum_{x \ne 0} \sup_{\eta} \left| c(\eta^x) - c(\eta) \right|,\\ 
\label{defep}
\epsilon &:=& \inf_{\eta} \left\{ c(\eta)+c(\eta^0) \right \},
\end{eqnarray}
where $\eta^x$ is the configuration obtained from $\eta$ by flipping the $x$-coordinate.

\subsubsection{Mixing for $\xi$}
\label{subsubsec:mixxi}

If $\xi$ is in the $M<\epsilon$ regime, then there is exponential decay of space-time correlations
(see Liggett~\cite{Li}, Section I.3). In fact, if $\xi$, $\xi'$ are two copies starting 
from initial configurations $\eta$, $\eta'$ and coupled according to the Vasershtein coupling, 
then, as was shown in Maes and Shlosman~\cite{MaSh}, the following estimate holds uniformly 
in $x \in \Z$ and in the initial configurations:
\begin{equation}
\label{expest1}
\mathbb{P}_{\eta, \eta'}\left( \xi_t(x) \ne {\xi'}_t(x)\right) 
\leq e^{-(\epsilon-M)t}.
\end{equation}
Since the system has uniformly bounded flip rates, it follows that there exist constants 
$K_1, K_2 \in (0,\infty)$, independent of $x \in \Z$ and of the initial configurations, such that
\begin{equation}
\label{expest2}
\mathbb{P}_{\eta, \eta'}\left( \exists \; s> t \text{ s.t. } \xi_s(x) 
\neq {\xi'}_s(x)\right) \le K_1 e^{-K_2t}.
\end{equation}
For $A \subset \Z \times \mathbb{R}_+$ measurable, let $\text{Discr}(A)$ be the event in 
which there is a discrepancy between $\xi$ and $\xi'$ in $A$, i.e., $\text{Discr}(A) 
:= \{\exists \; (x,t) \in A \colon\, \xi_t(x) \ne \xi'_t(x)\}$. Recall the definition of $C_R(m)$ 
in Section~\ref{subsec:notation}, and let $C_{R,t}(m) := C_R(m) \cap \Z \times [0,t]$. From (\ref{expest2}) we 
deduce that, for any fixed $m>0$ and $R \in \N_0$, there exist (possibly different) constants 
$K_1, K_2 \in (0,\infty)$ such that
\begin{equation}
\label{expestrestcone}
\mathbb{P}_{\eta,\eta'}(\text{Discr}(C_R(m) \setminus C_{R,t}(m))) \le K_1 e^{-K_2t}.
\end{equation}

\subsubsection{Mixing for $\Xi$}
\label{subsubsec:mixXi}

Bounds of the same type as \eqref{expest1}--\eqref{expestrestcone} hold for $\xi^{\pm}$, 
since $M = 0$ and $\epsilon > 0$ for independent spin-flips. Therefore, in order to have 
such bounds for the triple $\Xi$, we need only couple a pair $\Xi$, $\Xi'$ in such a 
way that each coordinate is coupled with its primed counterpart by the Vasershtein coupling. 
A set of coupling rates for $\Xi$, $\Xi'$ that accomplishes this goal is given 
in \eqref{ratesbigcoupling}, in Appendix \ref{appendix}.
Redefining $\text{Discr}(A):= \{\exists \; (x,t) \in A \colon\, \Xi_t(x) 
\neq \Xi'_t(x)\}$, by the previous results we see that (\ref{expestrestcone}) still 
holds for this coupling, with possibly different constants. As a consequence, we get the following lemma.

\begin{lemma}
\label{conecont}
Define $d(\eta,\eta'):= \sum_{x \in \Z}1_{\{\eta(x) \neq \eta'(x)\}}2^{-|x|-1}$. For 
any $m>0$ and $R \in \N_0$,
\begin{equation}
\label{eq:conecont}
\lim_{d(\Xi_0, \Xi'_0) \to 0}
\mathbb{P}_{\Xi_0,\Xi'_0}\big( \emph{Discr}(C_R(m)) \big)  = 0.
\end{equation}
\end{lemma}

\bpr
For any $t>0$, we may split $\text{Discr}(C_R(m)) = \text{Discr}(C_{R,t}(m)) \cup \text{Discr}
(C_R(m)\setminus C_{R,t}(m))$, so that
\begin{equation}
\label{estsplit}
\mathbb{P}_{\eta,\eta'}\big( \text{Discr}(C_R(m)) \big) 
\leq \mathbb{P}_{\eta,\eta'}\big( \text{Discr}(C_{R,t}(m)) 
+ \mathbb{P}_{\eta,\eta'}\big( \text{Discr}(C_R(m) \setminus C_{R,t}(m)) \big).
\end{equation}
Fix $\varepsilon > 0$. By (\ref{expestrestcone}), for $t$ large enough the second term in 
(\ref{estsplit}) is smaller than $\varepsilon$ uniformly in $\eta$, $\eta'$. For this 
fixed $t$, the first term goes to zero as $d(\eta,\eta') \to 0$, since $C_{R,t}(m)$ is 
contained in a finite space-time box and the coupling in (\ref{ratesbigcoupling}) is Feller 
with uniformly bounded total flip rates per site. (Note that the metric $d$ generates the 
product topology, under which the configuration space is compact.) Therefore $\limsup_{d(\eta,
\eta') \to 0} \mathbb{P}_{\eta, \eta '} \left( \text{Discr}(C_R(m)) \right) \leq 
\varepsilon$. Since $\varepsilon$ is arbitrary, (\ref{eq:conecont}) follows.
\epr

\subsubsection{Conditional mixing}
\label{subsubsec:condmix}

Next, we define an auxiliary process $\bar \Xi$ that, for each $L$, has the law of $\Xi$ 
conditioned on $\Gamma_L$ up to time $L$. We restrict to initial configurations $\eta \in 
\bar E$. In this case, $\bar \Xi$ is a process on $\left(\{0,1\}^{\Z \setminus \{0,1\}}\right)^3$ 
with rates that are equal to those of $\Xi$, evaluated with a trap at the origin. More precisely, 
for $\bar \eta \in \{0,1\}^{\Z\setminus\{0,1\}}$, denote by $(\bar \eta)_{1,0}$ the configuration 
in $\{0,1\}^\Z$ that is equal to $\bar \eta$ in $\Z\setminus\{0,1\}$ and has a trap at the origin. 
Then set $\bar C_x(\bar \eta) := C_x((\bar\eta)_{1,0})$, where $\bar C_x$ are the rates of 
$\bar \Xi$ and $C_x$ the rates of $\Xi$ at a site $x \in \Z$. Observe that the latter depend 
only on the middle configuration $\eta$, and not on $\eta^{\pm}$. These rates give the correct 
law for $\bar \Xi$ because $\Xi$ conditioned on $\Gamma_L$ is Markovian up to time $L$. Indeed,
the probability of $\Gamma_L$ does not depend on $\eta$ (for $\eta \in \bar E$) and, for $s < L$, 
$\Gamma_L = \Gamma_s \cap \theta_s \Gamma_{L-s}$. Thus, the rates follow by uniqueness. Observe 
that they are no longer translation-invariant.

Two copies of the process $\bar \Xi$ can be coupled analogously to $\Xi$ by restricting the rates in (\ref{ratesbigcoupling}) to $\bar E$. 
Since each coordinate of $\bar \Xi$ has similar properties as the 
corresponding coordinate in $\Xi$ (i.e., $\bar \xi^{\pm}$ are independent spin-flip systems and 
$\bar \xi$ is the in $M<\epsilon$ regime), it satisfies an estimate of the type
\begin{equation}
\label{expestauxproc}
\mathbb{\bar P}_{\eta, \eta'} \left( \text{Discr}([-t,t]\times \{t\}) \right) 
\leq K_1 e^{-K_2t} \quad \forall \;\; \eta, \eta' \in \bar E,
\end{equation}
for appropriate constants $K_1,K_2 \in (0,\infty)$. From this estimate we see that $d(\bar \Xi_t,\bar \Xi'_t) 
\to 0$ in probability as $t \to \infty$, uniformly in the initial configurations. By 
Lemma~\ref{conecont}, this is also true for $\mathbb{P}_{(\bar \Xi_t)_{1,0}, (\bar \Xi '_t)_{1,0}}
(\text{Discr}(C_R(m)))$. Since the latter is bounded, the convergence holds in $L_1$ as well, uniformly
in $\eta, \eta'$. 

\subsubsection{Proof of (H4)}
\label{proofH4ex1}

Let $f$ be a bounded function measurable in $\mathscr{C}_{R,\infty}(m)$ and estimate
\begin{align}
\label{h4mlessep}
&\left| \mathbb{E}_{\eta}\left[\theta_L f \mid \Gamma_L \right] 
- \mathbb{E}_{\eta'}\left[\theta_L f \mid \Gamma_L \right]\right| 
\leq 2 \|f\|_{\infty}\mathbb{P}_{\eta,\eta'}\big(\theta_L\text{Discr}(C_R(m))
\mid \Gamma_L \big) \nonumber \\
&\qquad \leq 2 \|f\|_{\infty} \sup_{\eta,\eta'}\mathbb{\bar E}_{\eta,\eta'}
\left[\mathbb{P}_{(\bar \Xi_L)_{1,0}, (\bar \Xi '_L)_{1,0}}
\left(\text{Discr}(C_R(m))\right)\right], 
\end{align}
where $\mathbb{\bar E}$ denotes expectation under the (coupled) law of $\bar \Xi$. Therefore (H4) 
follows with
\be{}
\Phi_L:= 2 \sup_{\eta,\eta'}\mathbb{\bar E}_{\eta,\eta'}
\left[\mathbb{P}_{(\bar \Xi_L)_{1,0}, (\bar \Xi '_L)_{1,0}}\left(\text{Discr}
(C_R(m))\right)\right],
\ee 
which converges to zero as $L \to \infty$ by the previous discussion. 
This is enough since $\kappa_L$ could be taken constant in the verification of (H3)(1), as we saw in \eqref{verifH31}.


\subsection{Example 2: subcritical dependence spread}
\label{sec:subcritdep}

In this section, we suppose that the rates $c(\eta)$ have a finite range of dependence
$r \in \N_0$. In this case, the system can be constructed via a graphical representation as 
follows.

\subsubsection{Graphical representation}
\label{subsubsec:graphrep}

For each $x \in \Z$, let $I^j_t(x)$ and $\Lambda^j_t(x)$ be independent Poisson processes 
with rates $c_j$ and $\lambda_j$ respectively, where $j=0,1$. At each event of 
$I^{j}_t(x)$, put a $j$-cross on the corresponding space-time point. At each event of 
$\Lambda^j(x)$, put two $j$-arrows pointing at $x$, one from each side, extending over the 
whole range of dependence. Start with an arbitrary initial configuration $\xi_0 \in \{0,1\}^{\Z}$. 
Then obtain the subsequent states $\xi_t(x)$ from $\xi_0$ and the Poisson processes by, at 
each $j$-cross, choosing the next state at site $x$ to be $j$ and, at at each $j$-arrow pair, 
choosing the next state to be $j$ if an independent Bernoulli($p_j(\theta_x\xi_s)$) trial 
succeeds, where $s$ is the time of the $j$-arrow event. This algorithm is well defined since, 
because of the finite range, up to each fixed positive time it can a.s.\ be performed locally.

Any collection of processes with the same range and with rates of the form (\ref{rates}) with 
$c_i$, $\lambda_i$ fixed ($i=0,1$) can be coupled together via this representation by fixing 
additionally for each site $x$ a sequence $(U_n(x))_{n \in \N}$ of independent Uniform$[0,1]$ 
random variables to evaluate the Bernoulli trials at $j$-arrow events. In particular, $\xi^{\pm}$ 
can be coupled together with $\xi$ in the graphical representation by noting that, for $\xi^-$, 
$p_0 \equiv 1$ and $p_1 \equiv 0$ and the opposite is true for $\xi^+$. For example, $\xi^-$ 
is the process obtained by ignoring all $1$-arrows and using all $0$-arrows. This gives the 
same coupling as the one given by the rates (\ref{ratescoupl}). In particular, we see that 
in this setting the events $\Gamma_L$ are given by (when $\xi_0 \in \bar E$)
\begin{equation}
\label{defgammalfiniterange}
\Gamma_L := \big\{I^0_L(0)=\Lambda^0_L(0)=I^1_L(1)=\Lambda^1_L(1)=0\big\}.
\end{equation}

\subsubsection{Coupling with a contact process}
\label{subsubsec:contproc}

We will couple $\Xi$ with a contact process $\zeta = (\zeta_t)_{t\geq 0}$ in the following way. 
We keep all Poisson events and start with a configuration $\zeta_0 \in \{i,h\}^\Z$, where $i$ 
stands for ``infected" and $h$ for ``healthy". We then interpret every cross as a recovery, 
and every arrow pair as infection transmission from any infected site within a neighborhood 
of radius $r$ to the site the arrows point to. This gives rise to a `threshold contact process' 
(TCP), i.e., a process with transitions at a site $x$ given by
\begin{equation}
\begin{array}{rcl}
i \rightarrow h & \text{ with rate } & c_0+c_1, \\ 
h \rightarrow i & \text{ with rate } & (\lambda_0 + \lambda_1) 
1_{\{ \exists \text{ infected site within range $r$ of $x$} \}}.
\end{array}
\end{equation}
In the graphical representation for $\xi$, we can interpret crosses as moments of memory loss 
and arrows as propagation of influence from the neighbors. Therefore, looking at the pair 
$(\Xi_t(x),\zeta_t(x))$, we can interpret the second coordinate being healthy as the first 
coordinate being independent of the initial configuration.

\begin{proposition}
\label{prop:couplingcp1}
Let $\underline i$ represent the configuration with all sites infected, and let $\Xi_0$, 
$\Xi'_0 \in E^3$. Couple $\Xi$, $\Xi'$ and $\zeta$ by fixing a realization of 
all crosses, arrows and uniform random variables, where $\Xi$ and $\Xi'$ are obtained from 
the respective initial configurations and $\zeta$ is started from $\underline i$. Then a.s.\
$\Xi_t(x)=\Xi'_t(x)$ for all $t > 0$ and $x \in \Z$ such that $\zeta_t(x) = h$.
\end{proposition}

\bpr
Fix $t>0$ and $x \in \Z$. With all Poisson and Uniform random variables fixed, an algorithm to find 
the state at $(x,t)$, simultaneously for any collection of systems of type (\ref{rates}) with 
fixed $c_i, \lambda_i$ and finite range $r$ from their respective initial configurations runs as follows. Find 
the first Poisson event before $t$ at site $x$. If it is a $j$-cross, then the state is 
$j$. If it is a $j$-arrow, then to decide the state we must evaluate $p_j$ and, therefore, 
we must first take note of the states at this time at each site within range $r$ of $x$, 
including $x$ itself. In order to do so, we restart the algorithm for each of these sites. 
This process ends when time $0$ or a cross is reached along every possible path from $(x,t)$ 
to $\Z\times\{0\}$ that uses arrows (transversed in the direction opposite to which they 
point) and vertical lines. In particular, if along each of these paths time $0$ is never 
reached, then the state at $(x,t)$ does not change when we change the initial configuration. 
On the other hand, time $0$ is not reached if and only if every path ends in a cross, which is 
exactly the description of the event $\{\zeta_t(x)=h\}$.
\epr

\subsubsection{Cone-mixing in the subcritical regime}
\label{subsubsec:subcritreg}

The process $(\zeta_t)_{t \geq 0}$ is stochastically dominated by a standard (linear) contact 
process (LCP) with the same range and rates. Therefore, if the LCP is subcritical, i.e., if 
$\lambda := (\lambda_0+\lambda_1)/(c_0+c_1) < \lambda_c$ where $\lambda_c$ is the critical 
parameter for the corresponding LCP, then the TCP will die out as well. 
Moreover, we have the following lemma:
\begin{lemma}
\label{lemma:subcritcp}
Let $A_t$ be the set of infected sites at time $t$. If $\lambda < \lambda_c$, 
then there exist positive constants $K_1, K_2, K_3, K_4$ such that
\begin{equation}
\mathbb{P}_{\underline i}\big(\exists\,s>t \colon\, A_s \cap [-K_1 e^{K_2s}, K_1 e^{K_2s}] \neq \emptyset \big) 
\leq K_3 e^{-K_4t}.
\end{equation}
\end{lemma}
\bpr This is a straightforward consequence of Proposition 1.1 in Aizenman-Jung \cite{AiJu07}.
\epr
According to Lemma~\ref{lemma:subcritcp}, the infection disappears exponentially fast around the 
origin. For $r=1$, a proof can be found in Liggett~\cite{Li}, Chapter VI, 
but it relies strongly on the nearest-neighbor nature of the interaction.

Let us now prove cone-mixing for $\xi$ when the rates are subcritical. 
Pick a cone $C_t$ with any inclination and tip at time $t$,
and let $\mathcal{H}_t := \{\text{all sites inside } C_t \text{ are healthy}\}$. This event is independent of $\xi_0$
and, because of Lemma~\ref{lemma:subcritcp}, has large probability if $t$ is large. 
Furthermore, by Proposition \ref{prop:couplingcp1}, on $\mathcal{H}_t$ the states of $\xi$ in $C_t$ are equal 
to a random variable that is independent $\xi_0$, which implies the cone-mixing property.

\subsubsection{Proof of (H4)}
\label{subsubsec:proofH4ex2}

In order to prove the \emph{conditional} cone-mixing property, we couple the 
conditioned process to a conditioned contact process as follows. First, let
\begin{equation}
\tilde \Gamma_L :=\left\{I_L^{j}(i)=\Lambda_L^{j}(i)=0 \colon\, j, i \in \{0,1\}\right\}.
\end{equation}

\begin{proposition}
\label{prop:couplingcp2}
Let $\hat i$ represent the configuration with all sites infected except for $\{0,1\}$, which are 
healthy. Let $\Xi_0$, $\Xi'_0 \in \bar{E}^3$. Couple $\Xi$, $\Xi'$ conditioned on 
$\Gamma_L$ and $\zeta$ conditioned on $\tilde \Gamma_L$ by fixing a realization of all crosses, 
arrows and uniform random variables as in Proposition~\ref{prop:couplingcp1} and starting, respectively, 
from $\Xi_0$, $\Xi'_0$ and $\hat i$, but, for $\Xi$ and $\Xi'$, remove the Poisson 
events that characterize $\Gamma_L$ and, for $\zeta$, remove all Poisson events up to time $L$ 
at sites $0$ and $1$, which characterizes $\tilde \Gamma_L$. Then a.s.\ $\Xi_t(x)=\Xi'_t(x)$ 
for all $t > 0$ and $x \in \Z$ such that $\zeta_t(x) = h$.
\end{proposition}

\bpr
On $\Gamma_L$, the states at sites $0$ and $1$ are fixed for time $[0,L]$. Therefore, in order 
to determine the state at $(x,t)$, we need not extend paths that touch $\{0,1\}\times[0,L]$: 
when every path from $(x,t)$ either ends in a cross or touches $\{0,1\}\times[0,L]$, the state 
at $(x,t)$ does not change when the initial configuration is changed in $\Z \setminus \{0,1\}$. 
But this is precisely the characterization of $\{\eta_t(x) = h\}$ on $\tilde \Gamma_L$ when 
started from $\hat i$.
\epr

The proof of (H4) is finished by noting that $(\eta_t)_{t\geq 0}$ starting from $\hat i$ and 
conditioned on $\tilde \Gamma_L$ is stochastically dominated by $(\eta_t)_{t\ge0}$ starting 
from $\underline i$. Therefore, by Lemma~\ref{lemma:subcritcp}, the ``dependence infection" 
still dies out exponentially fast, and we conclude as for the unconditioned cone-mixing.


\subsection{The sign of the speed}
\label{sec:signspeed}

For independent spin flips, we are able to characterize with the help of a coupling argument 
the regimes in which the speed is positive, zero or negative. By the stochastic domination 
described in Section~\ref{sec:sfsbr}, this gives us a criterion for positive (or negative) 
speed in the two classes addressed in Sections~\ref{sec:mlessep} and \ref{sec:subcritdep} above.

\subsubsection{Lipschitz property of the speed for independent spin-flip systems}

Let $\xi$ be an independent spin-flip system with rates $d_0$ and $d_1$ to flip to holes and 
particles, respectively. Since it fits both classes of IPS considered in Sections \ref{sec:mlessep} and \ref{sec:subcritdep}, by Theorem~\ref{mainthm} there exists a $w(d_0, d_1) \in \mathbb{R}$ 
that is the a.s.\ speed of the $(\infty,0)$-walk in this environment. This speed has the following
 local Lipschitz property.

\begin{lemma}
\label{signspeedisf}
Let $d_0, d_1, \delta > 0$. Then
\begin{equation}
w(d_0, d_1+\delta) - w(d_0, d_1 ) \geq \frac{d_0}{d_0 + d_1}\delta.
\end{equation}
\end{lemma}

\bpr
Our proof strategy is based on the proof of Theorem 2.24, Chapter VI in \cite{Li}.
Construct $\xi$ from a graphical representation by taking, 
for each site $x \in \Z$, two independent Poisson processes $N^i(x)$ with rates 
$d_i$, $i=0,1$, with each event of $N^i$ representing a flip to state $i$. For a fixed $\delta >0$, 
a second system $\xi^\delta$ with rates $d_0$ and $d_1 + \delta$ can be coupled to $\xi$ by starting 
from a configuration $\xi_0^\delta \ge \xi_0$ and adding to each site $x$ an independent Poisson process 
$N^\delta(x)$ with rate $\delta$, whose events also represent flips to particles, but only for 
$\xi^\delta$. Let us denote by $W$ and $W^\delta$ the walks in these respective environments. 
Under this coupling, $\xi \leq \xi^\delta$, so, by monotonicity, $W_t \le W_t^{\delta}$ 
for all $t \ge 0$ as well. We aim to prove that
\begin{equation}
\label{lips1}
\mathbb{E}_{\mu^{\delta}}\left[W^\delta_t\right] 
-\mathbb{E}_{\mu}\left[W_t\right] \geq \frac{d_0}{d_0+d_1}\delta t,
\end{equation}
where $\mu$ and $\mu^{\delta}$ are the equilibria of the respective systems. From this the 
conclusion will follow after dividing by $t$ and letting $t \to \infty$.

Define a third walk $W^*$ that is allowed to use one and only one event 
of $N^\delta$. More precisely, let $S$ be the first time when there is an event of $N^\delta$ at $W_S+1$. 
Take $W^*$ equal to $W$ on $[0,S)$ and, for times $\ge S$, let $W^*$ evolve by the same rules as $W$ but adding a particle at $W_S+1$ at time $S$, and using no more $N^\delta$ events. By construction, we have $W_t \le W_t^* \leq W_t^\delta$ $\forall$ $t \ge 0$.

Let $\eta_1 := \theta_{W_S}\xi_S \in \bar E$ and $\eta_2 := (\eta_1)^1$ be the configurations 
around $W_S$ and $W^*_{S-}$, respectively. Then
\begin{align}
\label{lips2}
\mathbb{E}_{\mu^{\delta}}\left[W^\delta_t\right] -\mathbb{E}_{\mu}\left[W_t\right] 
&\geq \mathbb{E}_{\mu}\left[W_t^*-W_t, S \le t\right] 
\geq \mathbb{E}_{\mu}\left[W_t^*-W_t, S \le t, \eta_1(2) =0\right] \nonumber \\ 
&= \mathbb{E}_{\mu}\left[\mathbb{E}_{\eta_1,\eta_2}
\left[ W^2_{t-S}-W^1_{t-S}\right], \eta_1(2)=0, S \le t \right],
\end{align}
where $W^i$, $i=1,2$ are copies of $W$ starting from $\eta_i$ and coupled via the graphical 
representation. We claim that, if  $\eta_1(2) = 0$, 
\begin{equation}\label{lips7}
\mathbb{E}_{\eta_1,\eta_2}\left[W^2_s - W^1_s \right] \ge 1 \quad \forall \; s \ge 0.
\end{equation}
Indeed, we will argue that the difference $W^2_s-W^1_s$ can only decrease 
when we flip all states of $\eta_1,\eta_2$ on $\Z_{\le -1}$ to particles 
and on $\Z_{\ge 2}$ to holes; but after 
doing these operations, we find that $W^2$ has the same distribution as $W^1+1$, which gives \eqref{lips7}.
It is enough to consider a single $x > 2$. Let $\tau := \inf\{t > 0 \colon\, N^0_t(x) + N^1_t(x) > 0 \} \wedge s$,
and put $T := \inf\{t > 0 \colon\, W^1_t = x-1\}$. There are two cases: either $T > \tau$ or not.
In the first case, $W^1_s$ remains constant if we set $\eta_{1,2}(x)=0$, while $W^2_s$ does not increase. 
In the second case, if $\eta_{1,2}(x)=0$, then $W^1_{T} = W^2_{T}$; 
but then they must remain equal thereafter since, for them to meet, the state at site $1$ must have flipped, 
and therefore they see the same configuration in the environment at time $T$. 
Hence, in this case, $W^2_s-W^1_s=0$ which is the minimum value, and our claim follows.

From \eqref{lips2} and \eqref{lips7} we get
\begin{equation}
\label{lips3}
\mathbb{E}_{\mu^{\delta}}\left[W^\delta_t\right] -\mathbb{E}_{\mu}\left[W_t\right] 
\geq \mathbb{P}_{\mu}\left(\eta_1(2)=0, S \le t \right).
\end{equation}
Consider the event $\{\eta_1(2)=0\} $. There are two possible situations: either at time $S$ 
the site $W_S +2$ was not yet visited by $W$, in which case $\eta_1(2)$ is still in equilibrium, 
or it was. In the latter case, let $s$ 
be the time of the last visit to this site before $S$. By geometrical constraints, at time 
$s$ only a hole could have been observed at this site, so the probability that its state at 
time $S$ is a hole is larger than at equilibrium, which is $d_0/(d_0+d_1)$. In other words,
\begin{equation}
\label{lips4}
\mathbb{P}_{\mu}\left(\eta_1(2)=0 \mid S, W_{[0,S]}\right) \geq \frac{d_0}{d_0+d_1},
\end{equation}
which, together with (\ref{lips3}) and the fact that $S$ has distribution Exp($\delta$), 
gives us
\begin{equation}
\label{lips5}
\mathbb{E}_{\mu^{\delta}}\left[W^\delta_t\right] -\mathbb{E}_{\mu}\left[W_t\right] 
\geq \frac{d_0}{d_0+d_1}\left(1-e^{\delta t}\right).
\end{equation}
Since $\delta$ is arbitrary, we may repeat the argument for systems with rates $d_1+ (k/n) 
\delta$, $n \in \N$ and $k =0,1,\ldots,n$, to obtain
\begin{align}
\label{lips5*}
\mathbb{E}_{\mu^{\delta}}\left[W^\delta_t\right] -\mathbb{E}_{\mu}\left[W_t\right] 
& \geq \frac{d_0}{d_0+d_1}n\left(1-e^{\delta t /n}\right),
\end{align}
and we get (\ref{lips1}) by letting $n \to \infty$.
\epr

\subsubsection{Sign of the speed}
\label{subsubsec:signspeed}

If $d_0 = d_1$, then $w = 0$, since by symmetry $W_t=-W_t$ in distribution. Hence we can summarize:

\begin{corollary}
\label{summarysignisf}
For an independent spin-flip system with rates $d_0$ and $d_1$,
\begin{equation}
\label{eqsumsignisf}
\begin{array}{ll}
w \ge \frac{d_0}{d_0+d_1}\left(d_1 - d_0\right) & \text{ if } d_1 > d_0,\\[0.2cm]
w = 0 & \text{ if } d_1 = d_0 ,\\[0.2cm]
w \le -\frac{d_1}{d_0+d_1}\left(d_0 - d_1\right) & \text{ if } d_1 < d_0.
\end{array}
\end{equation}
\end{corollary}

\noindent
Applying this result to the systems $\xi^{\pm}$ of Section~\ref{sec:sfsbr}, we obtain the following.

\begin{proposition}\label{critsignspeed}
Let $W$ be the random walk for the $(\infty,0)$-model in a spin-flip system 
with rates given by \eqref{rates}. Then, $\mathbb{P}_{\mu}$-a.s.,
\begin{equation}
\begin{array}{ll} 
\liminf_{t \to \infty} t^{-1}W_t \geq \frac{c_0+\lambda_0}{c_1+c_0+\lambda_0}
\left(c_1 - c_0 - \lambda_0 \right) & \text{ if } c_1 \ge c_0 + \lambda_0, \\[0.2cm]
\limsup_{t \to \infty} t^{-1}W_t \leq -\frac{c_1+\lambda_1}{c_0+c_1+\lambda_1}\left(c_0 - c_1 - \lambda_1 \right) 
& \text{ if } c_0 \ge c_1 + \lambda_1.
\end{array}
\end{equation}
\end{proposition}

This concludes the proof of Theorem~\ref{thm:cond} and the discussion of our two classes
of IPS's for the $(\infty,0)$-model. In Section~\ref{sec:oex} we give additional examples 
and discuss some limitations of our setting.


\section{Other examples}
\label{sec:oex}

We describe here three types of examples that satisfy our hypotheses: 
generalizations of the $(\alpha,\beta)$-model and of the $(\infty,0)$-model, and mixed models.
We also discuss an example that is beyond the reach of our setting.

\medskip\noindent
{\bf 1.\ Internal noise models.} 
For $x \in \Z\setminus\{0\}$ and $\eta \in E$, let $\pi_x(\eta)$ be functions with a finite range 
of dependence $R$. These are the rates to jump $x$ from the position of the walk. Let $\pi_x 
:= \sup_{\eta}\pi_x(\eta)$ and suppose that, for some $u>0$,
\begin{equation}
\label{expmomboundrates}
\sum_{x \in \Z\setminus\{0\}}e^{u|x|}\pi_x < \infty.
\end{equation}
This implies that also
\begin{equation}
\label{deftotrateinm}
\Pi := \sum_{x \in \Z\setminus\{0\}}\pi_x < \infty.
\end{equation}
The walk starts at the origin, and waits an independent Exponential($\Pi$) time $\tau$ until it 
jumps to $x$ with probability $\pi_x(\xi_\tau)/\Pi$. These probabilities do not necessarily sum 
up to one, so the walk may well stay at the origin. The subsequent jumps are obtained analogously, 
with $\xi_\tau$ substituted by the environment around the walk at the time of the attempted jump. 
It is clear that (A1)--(A3) hold. The walk has a bounded probability of standing still independently 
of the environment, and its jumps have an exponential tail. We take
\begin{equation}
\label{defgammainm}
\Gamma_L := \{\tau > L\}.
\end{equation}
By defining an auxiliary walk $(H_t)_{t \geq 0}$ that also tries to jump at time $\tau$, but only 
to sites $x>0$ with probability $\pi_x/\Pi$, we see that $W_t \leq H_t$ and that $H_t$ has 
properties analogous to the process defined in the proof of Lemma~\ref{sbisf}. Therefore, 
(H1)--(H3) are always satisfied for this model. Since $\Gamma_L$ is independent of $\xi$, 
(H4) is the (unconditional) cone-mixing property. Observe that $W_0=0$, so that $\bar \mu = \mu$. 
Therefore the LLN for this model holds in both examples discussed in Section \ref{sec:proofthm2}, and 
also for the IPS's for which cone-mixing was shown in Avena, den Hollander and Redig~\cite{AvdHoRe1}. 
The $(\alpha,\beta)$-model is an internal noise model with $R=0$ (the rates depend only on the 
state of the site where the walker is) and $\pi_x(\eta)=0$, except for $x = \pm1$, for which 
$\pi_{1}(1) = \alpha = \pi_{-1}(0)$ and $\pi_1(0)= \beta = \pi_{-1}(1)$.

\medskip\noindent
{\bf 2.\ Pattern models.} 
Take $\aleph$ to be a finite sequence of $0$'s and $1$'s, which we call a \emph{pattern}, and 
let $R$ be the length of this sequence. Take the environment $\xi$ to be of the same type used 
to define the $(\infty,0)$-walk. Let $q:\{0,1\}^R\setminus\{\aleph\} \to [0,1]$. The pattern 
walk is defined similarly as the $(\infty,0)$-walk, with the trap being substituted by the 
pattern, and a Bernoulli($q$) random variable being used to decide whether the walk jumps to 
the right or to the left. More precisely, let $\vartheta = (\xi_0(0),\ldots,\xi_0(R-1))$. 
If $\vartheta = \aleph$, then we set $W_0 = 0$, otherwise we sample $b_0$ as an independent Bernoulli($q(\vartheta$)) trial. 
If $b_0=1$, then $W_0$ is set to be the starting position 
of the first occurrence of $\aleph$ in $\xi_0$ to the right of the origin, while if $b_0=0$, 
then the first occurrence of $\aleph$ to the left of the origin is taken instead. Then the 
walk waits at this position until the configuration of one of the $R$ states to its right 
changes, at which time the procedure to find the jump is repeated with the environment as 
seen from $W_0$. Subsequent jumps are obtained analogously. The $(\infty,0)$-model is a 
pattern model with $\aleph := (1,0)$, $q(1,1) := 1$, $q(0,0):=0$ and $q(0,1):=1/2$.

For spin-flip systems given by (\ref{rates}), the pattern walk is defined and finite 
for all times, no matter what $\aleph$ is, the reasoning being exactly the same as for the 
$(\infty,0)$-walk. Also, it may be analogously defined so as to satisfy assumptions (A1)--(A3). Defining the events $\Gamma_L$ as
\begin{equation}
\label{deceventspattern}
\Gamma_L := \big\{\xi_s^{\pm}(j)=\xi_0^{\pm}(j) 
\; \forall \; s \in [0,L] \text{ and } j \in \{0,\ldots,R-1\}\big\},
\end{equation}
we may indeed, by completely analogous arguments, reobtain all the results of 
Section~\ref{sec:proofthm2}, so that hypotheses (H1)--(H4) hold and, therefore, 
the LLN as well.

\medskip\noindent
{\bf 3. Pattern models with extra jumps.}
Examples of models that fall into our setting and for which the events $\Gamma_L$ depend
non-trivially both on $\xi$ and $Y$ can be constructed by taking a pattern model and adding noise 
in the form of non-zero jump rates while sitting on the pattern.
More precisely, add to $Y$ an independent Poisson process $N$ with positive rate
and let $W$ jump also at events of $N$ but with the same jump mechanism, i.e.,
choosing the sign of the jump according to the result of a Bernoulli($q$) random variable, and the displacement using the pattern.
Taking $\Gamma_L := \Gamma^{\aleph}_L \cap \{N_L=0\}$, where $\Gamma^{\aleph}_L$ is the corresponding event for the pattern model, 
we may check that, for the two examples of dynamic random environments 
considered in Theorem~\ref{thm:cond}, (A1)--(A3) and (H1)--(H4) are all verified.

\medskip\noindent
{\bf 4. Mixtures of pattern and internal noise models.}
Another class of models with nontrivial dependence structure for the renegeration-inducing events can be constructed as follows.
Let $W^0$ be an internal noise model and $W^1$ a pattern model (with or without extra jumps) on the the same random environment $\xi$
and let $Y^i$, $i \in \{0,1\}$, be the corresponding random elements associated to each model.
Let $X = (X)_{n \in \N}$ be a sequence of i.i.d.\ Bernoulli($p$) random variables independent of all the rest, where $p \in (0,1)$. 
Then the mixture is the model for which the dynamics associated to $i \in \{0,1\}$ are applied in the time interval $[n-1,n)$ when 
$X_n = i$. Note that this model will have deterministic jumps.

Letting $Y := \left(Y^0,Y^1, X \right)$ where $Y^i$ is the corresponding
random element associated to the model $i$, 
it is easily checked that this model indeed falls into our setting.

Choosing
\begin{equation}
\label{defgammamix}
\Gamma_L := \Gamma_L^1 \cap \{X_k = 1, k = 1,\ldots,L\}
\end{equation}
where $\Gamma^1_L$ is the corresponding event for the pattern model, 
it is not hard to verify, using the results of Section~\ref{sec:proofthm2}, that, for the two classes of random environments considered in
Theorem~\ref{thm:cond}, the mixed model satisfies (A1)--(A3) and (H1)--(H4).

\medskip\noindent
{\bf An open example.}
We will close with an example of a model that does not satisfy the hypotheses of our 
LLN (in dynamic random environments given by spin-flip systems). When $\xi(0)=j$, let 
$C^j$ be the cluster of $j$'s around the origin. Define jump rates for the walk as follows:
\begin{equation}
\label{ratescounterex}
\begin{array}{lcl}
\pi_1(\eta) & = & \left\{\begin{array}{lll} 
               |C^1| & \text{ if } & \eta(0)=1,\\
               |C^0|^{-1} & \text{ if } & \eta(0)=0,
               \end{array}\right.\\
\pi_{-1}(\eta) & = & \left\{ \begin{array}{lll} 
               |C^0| & \text{ if } & \eta(0)=0,\\
               |C^1|^{-1} & \text{ if } & \eta(0)=1.
               \end{array} \right.
\end{array}
\end{equation}
Even though this looks like a fairly natural model, 
it does not satisfy (A2). It also won't satisfy (H1) and (H2) together for any reasonable random environment, 
which is actually the hardest obstacle.
The problem is that, while we are able to transport 
a.s.\ properties of the equilibrium measure to the measure of the environment as seen from the walk, 
we cannot control the distortion in events 
of positive measure. Thus, even if $\Gamma_L$ has positive probability at time zero, there 
is no a priori guarantee that it will have an appreciable probability 
from the point of view of the walk at later times. 
Because of this, we cannot implement our regeneration strategy, and our proof of the LLN breaks down.

\newpage
\appendix
\section{Appendix: coupling rates}
\label{appendix}

Here we give the rates for a coupling between $\Xi$ and $\Xi'$, mentioned in Section \ref{subsubsec:mixXi}, 
such that corresponding pairs of coordinates are distributed according to the Vasershtein coupling. 
Let $\eta$, $\eta'$ be the state of the middle coordinates $\xi$ and $\xi'$; the states outside the origin of the other coordinates play no role. Then the flip rates at the origin are given schematically by
\begin{equation}
\label{ratesbigcoupling}
\begin{array}{ccl}
(000)(000) & \to & \left\{ \begin{array}{ll}
                   (111)(111) & c_1, \\
                   (011)(011) & c(\eta) \wedge c(\eta') - c_1,\\
                   (011)(001) & c(\eta) - c(\eta) \wedge c(\eta'),\\
                   (001)(011) & c(\eta') - c(\eta) \wedge c(\eta'),\\
                   (001)(001) & c_1+\lambda_1 - c(\eta) \vee c(\eta'),
                   \end{array}\right.\\
(001)(001) & \to & \left\{ \begin{array}{ll}
                   (111)(111) & c_1, \\
                   (011)(011) & c(\eta) \wedge c(\eta') - c_1,\\
                   (011)(001) & c(\eta)-c(\eta) \wedge c(\eta'),\\
                   (001)(011) & c(\eta')-c(\eta) \wedge c(\eta'),\\
                   (000)(000) & c_0,
                   \end{array}\right.\\
(001)(011) & \to & \left\{ \begin{array}{ll}
                   (111)(111) & c_1, \\
                   (011)(011) & c(\eta) - c_1,\\
                   (001)(001) & c(\eta') - c_0,\\
                   (000)(000) & c_0,
                   \end{array}\right.\\
(000)(001) & \to & \left\{ \begin{array}{ll}
                   (111)(111) & c_1, \\
                   (011)(011) & c(\eta) \wedge c(\eta') - c_1,\\
                   (011)(001) & c(\eta)-c(\eta) \wedge c(\eta'),\\
                   (001)(011) & c(\eta')-c(\eta) \wedge c(\eta'),\\
                   (001)(001) & c_1+\lambda_1 - c(\eta) \vee c(\eta'),\\
                   (000)(000) & c_0,
                   \end{array}\right.\\
(000)(011) & \to & \left\{ \begin{array}{ll}
                   (111)(111) & c_1, \\
                   (011)(011) & c(\eta) - c_1,\\
                   (001)(011) & c_1+\lambda_1-c(\eta),\\
                   (000)(000) & c_0,\\
                   (000)(001) & c(\eta')-c_0,
                   \end{array}\right.\\
(000)(111) & \to & \left\{ \begin{array}{ll}
                   (111)(111) & c_1, \\
                   (011)(111) & c(\eta) - c_1,\\
                   (001)(111) & c_1+\lambda_1-c(\eta),\\
                   (000)(000) & c_0,\\
                   (000)(001) & c(\eta')-c_0,\\
                   (000)(011) & c_0 + \lambda_0 - c(\eta').
                   \end{array}\right.
\end{array}
\end{equation}
The other transitions, starting from 
\be{}
(111)(111), \quad (011)(011), \quad (011)(001) , \quad (111)(011), \quad (111)(001) \; \text{ and } \; (111)(000),
\ee
can be obtained from the ones in \eqref{ratesbigcoupling} by symmetry, 
by exchanging the roles of $\eta$/$\eta'$ or of particles/holes.


\end{document}